\documentclass[aop,seceqn,MSNbibl,citesort,dvips]{arximspdf}
\usepackage{mathbh}

\doi{10.1214/10-AOP556}
\volume{39}
\issue{2}
\pubyear{2011}
\firstpage{471}
\lastpage{506}

\makeatletter

\newtheorem{theorem}{Theorem}[section]
\newtheorem{corollary}[theorem]{Corollary}
\newproclaim{definition}[theorem]{Definition}
\newtheorem{lemma}[theorem]{Lemma}
\newproclaim{remark}[theorem]{Remark}

\newcommand{\one}{\mathbh{1}}
\newcommand{\dd}{{d}}
\newcommand{\ee}{{e}}
\newcommand{\eq}{\mathit{eq}}

\makeatother

\begin{document}
\begin{frontmatter}

\title{Harmonic functions, \textup{h}-transform and large deviations for random
walks in random environments in dimensions four and higher}
\runtitle{Large deviations for RWRE}

\begin{aug}
\author[A]{\fnms{Atilla} \snm{Yilmaz}\corref{}\ead[label=e1]{atilla@math.berkeley.edu}}
\runauthor{A. Yilmaz}
\affiliation{University of California, Berkeley}
\address[A]{Department of Mathematics\\
University of California, Berkeley\\
Berkeley, California 94720-3840\\
USA\\
\printead{e1}} 
\end{aug}

\received{\smonth{1} \syear{2010}}
\revised{\smonth{4} \syear{2010}}

%
\begin{abstract}
We consider large deviations for nearest-neighbor random walk in a
uniformly elliptic i.i.d. environment on $\mathbb{Z}^d$. There exist
variational formulae for the quenched and averaged rate functions
$I_q$ and $I_a$, obtained by Rosenbluth and Varadhan, respectively.
$I_q$ and $I_a$ are not identically equal. However, when $d\geq4$ and
the walk satisfies the so-called (\textbf{T}) condition of Sznitman,
they have been previously shown to be equal on an open set
$\mathcal{A}_{\eq}$.

For every $\xi\in\mathcal{A}_{\eq}$, we prove the existence of a
positive solution to a Laplace-like equation involving $\xi$ and the
original transition kernel of the walk. We then use this solution to
define a new transition kernel via the h-transform technique of Doob.
This new kernel corresponds to the unique minimizer of Varadhan's
variational formula at $\xi$. It also corresponds to the unique
minimizer of Rosenbluth's variational formula, provided that the latter
is slightly modified.
\end{abstract}

%
\begin{keyword}[class=AMS]
\kwd{60K37}
\kwd{60F10}
\kwd{82C41}.
\end{keyword}
\begin{keyword}
\kwd{Random walk}
\kwd{random environment}
\kwd{large deviations}
\kwd{harmonic functions}
\kwd{Doob h-transform}.
\end{keyword}

\end{frontmatter}

\section{Introduction}

\subsection{The model}

Let $(e_i)_{i=1}^d$ be the canonical basis for the $d$-dimensional
integer lattice $\mathbb{Z}^d$ with $d\geq1$. Consider a discrete-time
Markov chain on $\mathbb{Z}^d$ with nearest-neighbor steps, that is,
with steps in $U:=\{\pm e_i\}_{i=1}^d$. For every $x\in\mathbb{Z}^d$
and $z\in U$, denote the transition probability from $x$ to $x+z$ by
$\pi(x,x+z)$ and refer to the transition vector $\omega_x:=(\pi
(x,x+z))_{z\in U}$ as the \textit{environment} at $x$. If the
environment $\omega:=(\omega_x)_{x\in\mathbb{Z}^d}$ is sampled from a
probability space $(\Omega,\mathcal{B},\mathbb{P})$, then this process
is called \textit{random walk in a random environment} (RWRE). Here,
$\mathcal{B}$ is the Borel $\sigma$-algebra corresponding to the
product topology.

The environment is said to be \textit{uniformly elliptic} if
%
%
\begin{equation}\label{ellipticity}\qquad
\mbox{there exists a $\delta>0$ such that $\pi(0,z)\geq\delta$ for
every $\omega\in\Omega$ and $z\in U$.}
\end{equation}

For every $y\in\mathbb{Z}^d$, define the shift $T_y$ on $\Omega$ by
$ (T_y\omega)_x:=\omega_{x+y}$. Throughout this paper, we
will assume that $\mathbb{P}$ is stationary and ergodic under $
(T_z )_{z\in U}$. This condition is clearly satisfied when
%
%
\begin{equation}\label{i.i.d.}
\omega=(\omega_x)_{x\in\mathbb{Z}^d}\qquad\mbox{is an i.i.d. collection.}
\end{equation}

For every $x\in\mathbb{Z}^d$ and $\omega\in\Omega$, the Markov chain
with environment $\omega$ induces a probability measure $P_x^\omega$ on
the space of paths starting at $x$. Statements about $P_x^\omega$ that
hold for $\mathbb{P}$-a.e. $\omega$ are referred to as \textit
{quenched}. Statements about the semidirect product $P_x:=\mathbb
{P}\otimes P_x^\omega$ are referred to as \textit{averaged} (or
\textit{annealed}). Expectations under $\mathbb{P}, P_x^\omega$ and
$P_x$ are
denoted by $\mathbb{E}, E_x^\omega$ and $E_x$, respectively.

See \cite{ZeitouniSurvey06} for a survey of results on RWRE.

Because of the extra layer of randomness in the model, the standard
questions of recurrence versus transience, the law of large numbers
(LLN), the central limit theorem (CLT) and the large deviation
principle (LDP)---which have well-known answers for classical random
walk---become hard. However, it is possible, by taking the
\textit{point of view of the particle}, to treat the two layers of randomness
as one: if we denote the random path of the particle by $X:=(X_n)_{n\geq
0}$, then $(T_{X_n}\omega)_{n\geq0}$ is a Markov chain (referred to as
the \textit{environment Markov chain}) on $\Omega$ with transition
kernel $\overline{\pi}$ given by
\[
\overline{\pi}(\omega,\omega'):=\sum_{z\dvtx T_z\omega=\omega'}\pi(0,z).
\]
This is a standard approach in the study of random media; see, for
example, \cite{KipnisVaradhan,Kozlov} or \cite{PapaVaradhan}.

Instead of viewing the environment Markov chain as an auxiliary
construction, one can introduce it first and then deduce the particle
dynamics from it.
\begin{definition}\label{ortamkeli}
A function $\hat{\pi}\dvtx\Omega\times U\to\mathbb{R}^+$ is said to be an
``environment kernel'' if $\hat{\pi}(\cdot,z)$ is $\mathcal
{B}$-measurable for each $z\in U$ and $\sum_{z\in U}\hat{\pi}(\cdot,z)=1$.
It can be viewed as a transition kernel on $\Omega$ via the following
identification:
\[
\overline{\pi}(\omega,\omega'):=\sum_{z\dvtx T_z\omega=\omega'}\hat{\pi
}(\omega,z).
\]
Given $x\in\mathbb{Z}^d$, $\omega\in\Omega$ and any environment kernel
$\hat{\pi}$, the quenched probability measure $P_x^{\hat{\pi},\omega}$
on the space of particle paths $(X_n)_{n\geq0}$ starting at $x$ in
environment $\omega$ is defined by setting $P_x^{\hat{\pi},\omega}
(X_o=x )=1$ and
\[
P_x^{\hat{\pi},\omega} (X_{n+1}=y+z |X_n=y )=\hat
{\pi}(T_y\omega,z)
\]
for all $n\geq0$, $y\in\mathbb{Z}^d$ and $z\in U$. The semidirect
product $P_x^{\hat{\pi}}:=\mathbb{P}\otimes P_x^{\hat{\pi},\omega}$ is
referred to as the averaged measure and expectations under $P_x^{\hat
{\pi},\omega}$ and $P_x^{\hat{\pi}}$ are denoted by $E_x^{\hat{\pi
},\omega}$ and $E_x^{\hat{\pi}}$, respectively.
\end{definition}
%

\subsection{Summary of results}

In this paper, we will focus on the large deviation properties of
multidimensional RWRE. Section \ref{kurdungel} is a detailed survey of
the previous results on this topic that are relevant to our purposes.
The precise statements of our results are postponed to Section \ref
{inisoyler} because they rely heavily on the notation and theorems
given in Section \ref{kurdungel}.

In this subsection, we will provide a short and less technical
description of the key theorems in Section \ref{kurdungel}. References
will be omitted for the sake of brevity. We will then highlight our
main results.

%

\subsubsection{Summary of previous results}

In the case of quenched RWRE, the LDP holds for the mean velocity
$X_n/n$ of the particle. 
Rosenbluth gives a variational formula for the corresponding rate
function $I_q$. 
For any $\xi\in\mathbb{R}^d$, $I_q(\xi)$ is equal to the infimum of
$H(\hat{\pi},\mathbb{Q})$, where $H(\cdot)$ is a relative entropy and
$(\hat{\pi},\mathbb{Q})$ varies over all pairs such that: (i) $\hat{\pi
}$ is an environment kernel; (ii) $\mathbb{Q}$ is a $\hat{\pi
}$-invariant probability measure on $\Omega$; (iii) $\mathbb{Q}\ll
\mathbb{P}$ on $\mathcal{B}$; (iv) the asymptotic mean velocity of the
walk induced by $(\hat{\pi},\mathbb{Q})$ is equal to $\xi$.

For averaged walks in i.i.d. environments, Varadhan proves the LDP for
$X_n/n$ 
and gives yet another variational formula for the corresponding rate
function $I_a$. 
For any $\xi\neq0$, $I_a(\xi)$ is the infimum of $\mathfrak{I}_a(\alpha
)$, where $\mathfrak{I}_a(\cdot)$ is a relative entropy [not equal to
$H(\cdot)$] and $\alpha$ varies over all $\mathbb{Z}^d$-valued
transient processes with stationary and ergodic increments in $U$ such
that the mean drift of $\alpha$ is equal to $\xi$.

It is easily shown that (i) $I_a\leq I_q$ and (ii) $I_q$, $I_a$ are not
identically equal. When $d\geq4$ and the walk satisfies the so-called
(\textbf{T}) condition of Sznitman, $I_q$ and $I_a$ are known to be
strictly convex, analytic and equal on an open set $\mathcal{A}_{\eq}$.
At every $\xi\in\mathcal{A}_{\eq}$, Varadhan's variational formula for
$I_a(\xi)$ has a unique minimizer.

%
\subsubsection{Summary of our results}

We will assume that the environment is i.i.d., $d\geq4$ and the
(\textbf{T}) condition of Sznitman holds. For every $\xi\in\mathcal{A}_{\eq}$,
we will prove the existence of an $h(\theta,\cdot)\in L^2(\mathbb{P})$
that solves a certain equation involving $\theta:=\nabla I_a(\xi)$ and
the original kernel $\pi$ of the walk; see (\ref{carbonic}). Since
(\ref{carbonic}) resembles the Laplace equation, we will refer to
$h(\theta,\cdot)$ as \textit{harmonic}. We will then use $h(\theta,\cdot
)$ to define a new environment kernel $\hat{\pi}^\theta$ via the
h-transform technique of Doob; see~(\ref{azdinyn}).

For every $\xi\in\mathcal{A}_{\eq}$, we will prove the existence of a
probability measure $\mathbb{Q}_\xi$ on $\Omega$ that is $\hat{\pi
}^\theta$-invariant. The pair $(\hat{\pi}^\theta,\mathbb{Q}_\xi)$
corresponds to a stationary Markov chain with values in $\Omega$. This
Markov chain induces a $\mathbb{Z}^d$-valued transient process $\mu_\xi
^\infty$ with stationary and ergodic increments in $U$. We will show
that $\mu_\xi^\infty$ is the unique minimizer of Varadhan's variational
formula for $I_a(\xi)$.

The pair $(\hat{\pi}^\theta,\mathbb{Q}_\xi)$ is a natural minimizer
candidate for Rosenbluth's variational formula for $I_q(\xi)$. However,
it is not known whether $\mathbb{Q}_\xi\ll\mathbb{P}$ on $\mathcal{B}$.
We will resolve this issue by slightly modifying Rosenbluth's formula
so that the infimum of $H(\cdot)$ will be taken over a larger class of
pairs. Finally, we will show that $(\hat{\pi}^\theta,\mathbb{Q}_\xi)$
is the unique minimizer of this new formula.

%
\section{Previous results on large deviations for RWRE}\label{kurdungel}

\subsection{The quenched LDP}

Recall that a sequence $ (Q_n )_{n\geq1}$ of probability
measures on a topological space $\mathbb{X}$ satisfies\vadjust{\goodbreak} the
\textit{large deviation principle} (LDP) with rate function $I\dvtx
\mathbb{X}\to
[0,\infty]$ if $I$ is lower semicontinuous and, for any measurable set
$G$,
\[
-\inf_{x\in G^o}I(x)\leq\liminf_{n\to\infty}\frac{1}{n}\log Q_n(G)\leq
\limsup_{n\to\infty}\frac{1}{n}\log Q_n(G)\leq-\inf_{x\in\overline{G}}I(x).
\]
Here, $G^o$ is the interior of $G$ and $\overline{G}$ its closure. See
\cite{DemboZeitouniBook} for general background regarding large deviations.

In this paper, the following theorem will be referred to as the
quenched (level-1) LDP.
\begin{theorem}[(Quenched LDP)]\label{qLDPgeneric}
Assume (\ref{ellipticity}). For $\mathbb{P}$-a.e. $\omega$, $
(P_o^\omega(\frac{X_n}{n}\in\cdot ) )_{n\geq1}$
satisfies the LDP with a deterministic and convex rate function $I_q$.
(The subscript stands for ``quenched.'')
\end{theorem}

Greven and den Hollander \cite{GdH94} prove Theorem \ref{qLDPgeneric}
for walks on $\mathbb{Z}$ in i.i.d. environments. They provide a
formula for $I_q$ and show that its graph typically has flat pieces.
Comets, Gantert and Zeitouni \cite{CGZ00} generalize the results in
\cite{GdH94} to stationary and ergodic environments.

For $d\geq1$, the first result on quenched large deviations is given by
Zerner \cite{Zerner98}. He uses a subadditivity argument for certain
passage times to prove Theorem \ref{qLDPgeneric} in the case of
\textit{nestling} walks in i.i.d. environments.
\begin{definition}\label{nonnennon}
RWRE is said to be \textit{nonnestling} relative to a unit vector $\hat
{u}\in\mathcal{S}^{d-1}$ if
%
%
\begin{equation}\label{nonmumu}
\operatorname{ess}\inf_{\mathbb{P}}\sum_{z\in U}\pi(0,z)\langle z,\hat
{u}\rangle>0.
\end{equation}
It is said to be \textit{nestling} if it is not nonnestling relative to
any unit vector. In the latter case, the convex hull of the support of
the law of $\sum_{z}\pi(0,z)z$ contains the origin.
\end{definition}

By a more direct use of the subadditive ergodic theorem, Varadhan
\cite{Varadhan03} drops the nestling assumption and generalizes Zerner's
result to stationary and ergodic environments. The drawback of these
approaches is that they do not lead to any formula for the rate function.

Kosygina, Rezakhanlou and Varadhan \cite{KRV06} consider diffusions on
$\mathbb{R}^d$ (with \mbox{$d\geq1$}) in stationary and ergodic environments.
They prove the analog of Theorem~\ref{qLDPgeneric} via a minimax
argument and provide a variational formula for the quenched rate
function. Rosenbluth \cite{RosenbluthThesis} adapts their work to the
context of RWRE. [See (\ref{level1ratetilde}) below for Rosenbluth's
variational formula for $I_q$.]
%

\subsection{The quenched level-2 LDP and Rosenbluth's variational formula}

The minimax argument of Kosygina et al. \cite{KRV06} can be
generalized to establish a quenched LDP for the so-called \textit{pair
empirical measure} of the environment Markov chain. Below, we introduce
some notation in order to give the precise statement of this theorem.

For any measurable space $(Y,\mathcal{F})$, write $M_1(Y,\mathcal{F})$
[or simply $M_1(Y)$ whenever no confusion occurs] for the space of
probability measures on $(Y,\mathcal{F})$. Consider the random walk
$X=(X_n)_{n\geq0}$ on $\mathbb{Z}^d$ in a stationary and ergodic
environment, let $Z_n=X_n-X_{n-1}$ and focus on
\[
\nu_{n,X} := \frac{1}{n}\sum_{k=0}^{n-1}\one_{T_{X_k}\omega,Z_{k+1}},
\]
which is a random element of $M_1(\Omega\times U)$. The map $(\omega
,z)\mapsto(\omega,T_z\omega)$ embeds $M_1(\Omega\times U)$ into
$M_1(\Omega\times\Omega)$ and we therefore refer to $\nu_{n,X}$ as the
pair empirical measure of the environment Markov chain. For any $\mu\in
M_1(\Omega\times U)$, define the probability measures $(\mu)^1$ and
$(\mu)^2$ on $\Omega$ by
%
%
\begin{equation}\label{albaytey}
\dd(\mu)^1(\omega):=\sum_{z\in U}\dd\mu(\omega,z) \quad\mbox
{and}\quad \dd(\mu)^2(\omega):=\sum_{z\in U}\dd\mu
(T_{-z}\omega,z),
\end{equation}
respectively, which are the marginals of $\mu$ when $\mu$ is seen as an
element of $M_1(\Omega\times\Omega)$. With this notation, let
\begin{eqnarray*}
M_1'(\Omega\times U)&:=& \biggl\{\mu\in M_1(\Omega\times U)\dvtx(\mu)^1=(\mu
)^2\ll\mathbb{P},\\
&&\hspace*{33.1pt} \frac{\dd\mu(\cdot,z)}{\dd(\mu)^1(\cdot
)}>0\mbox{ for every }z\in U \biggr\}.
\end{eqnarray*}
\begin{theorem}[(Quenched level-2 LDP, Yilmaz
\cite{YilmazQuenched})]\label{level2LDP} Assume (\ref{ellipticity}).
For $\mathbb{P}$-a.e. $\omega$, $(P_o^\omega(\nu_{n,X}\in\cdot))_{n\geq
1}$ satisfies the LDP with the rate function $\mathfrak{I}_q^{**}$, the
double convex conjugate of $\mathfrak{I}_q\dvtx M_1(\Omega\times U)\to
\mathbb
{R}$ given by
%
%
\begin{equation}\label{level2ratetilde}
\mathfrak{I}_q(\mu)= \cases{
\displaystyle\int_{\Omega}\sum_{z\in U} \dd\mu(\omega,z)\log\frac{\dd\mu(\omega
,z)}{\dd(\mu)^1(\omega)\pi(0,z)},&\quad if $\mu\in
M_1'(\Omega\times U)$,\cr
\infty, &\quad otherwise.}\hspace*{-32pt}
\end{equation}
\end{theorem}

Rosenbluth's quenched LDP result is a corollary of Theorem
\ref{level2LDP}. Indeed, for any $\mu\in M_1(\Omega\times U)$, set
%
%
\begin{equation}\label{ximu}
\xi_{\mu}:=\int\sum_{z\in U}\dd\mu(\omega,z)z.
\end{equation}
For any $\xi\in\mathbb{R}^d$, define
%
%
\begin{equation}\label{Axi}
A_\xi:=\{\mu\in M_1(\Omega\times U)\dvtx\xi_{\mu}=\xi\}.
\end{equation}
With this notation,
%
%
\begin{eqnarray}\label{level1rate}
I_q(\xi)&=&\inf_{\mu\in A_\xi} \mathfrak{I}_q^{**}(\mu)\\
\label{level1ratetilde}
&=&\inf_{\mu\in A_\xi} \mathfrak{I}_q(\mu).
\end{eqnarray}
Here, (\ref{level1rate}) follows from Theorem \ref{level2LDP} via the
so-called contraction principle (see~\cite{DemboZeitouniBook}).
Note that, even though $\mathfrak{I}_q$ is convex, it may not be lower
semicontinuous (see Appendix A of \cite{YilmazQuenched} for an
example). Therefore, $\mathfrak{I}_q^{**}$ is not equal to $\mathfrak
{I}_q$ in general. Nevertheless, (\ref{level1ratetilde}) is valid (see
\cite{YilmazQuenched}) and it is precisely equal to the variational
formula obtained by Rosenbluth in \cite{RosenbluthThesis}.

%

\subsection{The quenched level-3 LDP}

Theorem \ref{level2LDP} can be generalized to establish a quenched LDP
for the \textit{empirical process}
\[
\nu_{n,X}^\infty:= \frac{1}{n}\sum_{k=0}^{n-1}\one_{T_{X_k}\omega
,Z_{k+1}^\infty},
\]
which is a random element of $M_1(\Omega\times U^\mathbb{N})$. Here,
$Z_{k+1}^\infty$ is shorthand notation for $(Z_{k+i})_{i\geq1}$.
\begin{theorem}[(Quenched level-3 LDP, Rassoul-Agha and Sepp\"al\"ainen
\cite{FirasTimoLevel3})]\label{level3LDP} Assume (\ref{ellipticity}).
For $\mathbb{P}$-a.e. $\omega$, $(P_o^\omega(\nu_{n,X}^\infty\in\cdot
))_{n\geq1}$ satisfies the LDP with a deterministic and convex rate function
$I_{q,3}\dvtx M_1(\Omega\times U^\mathbb{N})\to\mathbb{R}$.
\end{theorem}

Rassoul-Agha and Sepp\"al\"ainen actually obtain this result in greater
generality, namely for bounded step size walks satisfying a weak
ellipticity condition (see~\cite{FirasTimoLevel3}). Also, they show
that, just as in Theorem \ref{level2LDP}, the rate function $I_{q,3}$
is the lower semicontinuous regularization of a relative entropy. We
choose not to state the precise formula of $I_{q,3}$ here, partly in
order to keep the notation simple and partly because we will not need
it in what follows.

%

\subsection{The averaged LDP and Varadhan's variational formula}

In this paper, the following theorem will be referred to as the
averaged (level-1) LDP.

\begin{theorem}[(Averaged LDP)]\label{aLDPgeneric}
Assume (\ref{ellipticity}) and (\ref{i.i.d.}). $ (P_o (\frac
{X_n}{n}\in\cdot ) )_{n\geq1}$ satisfies the LDP with a
convex rate function $I_a$ (the subscript stands for ``averaged'').
\end{theorem}

Comets et al. \cite{CGZ00} prove Theorem \ref{aLDPgeneric} for $d=1$
and obtain the following variational formula for $I_a$:
%
%
\begin{equation}\label{ahubay}
I_a(\xi)=\inf_{\mathbb{Q}} \{I_q^\mathbb{Q}(\xi) + |\xi|h_s
(\mathbb{Q} |\mathbb{P} ) \}.
\end{equation}
Here, the infimum is over all stationary and ergodic probability
measures on $\Omega$, $I_q^\mathbb{Q}(\cdot)$ denotes the rate function
for the quenched LDP when the environment measure is $\mathbb{Q}$ and
$h_s (\cdot|\cdot)$ is specific relative entropy.
Similarly to the quenched picture, the graph of $I_a$ is shown
typically to have flat pieces.

Varadhan \cite{Varadhan03} proves Theorem \ref{aLDPgeneric} for any
$d\geq1$. He gives yet another variational formula for $I_a$. Below, we
introduce some notation in order to write down this formula.

An infinite path $ (x_i )_{i\leq0}$ with nearest-neighbor
steps $x_{i+1}-x_i$ is said to be in $W_\infty^{\mathrm{tr}}$ if
$x_o=0$ and $\lim_{i\to-\infty}|x_i|=\infty$. For any $w\in W_\infty
^{\mathrm{tr}}$, let $n_o$ be the number of times $w$ visits the
origin, excluding the last visit. By the transience assumption, $n_o$
is finite. For any $z\in U$, let $n_{o,z}$ be the number of times $w$
jumps to $z$ after a visit to the origin. Clearly, $\sum_{z\in
U}n_{o,z}=n_o$. If the averaged walk starts from time $-\infty$ and its
path $ (X_i )_{i\leq0}$ up to the present is conditioned to be
equal to $w$, then the probability of the next step being equal to $z$ is
%
%
\begin{equation}\label{ozgurevren}
q(w,z):=\frac{\mathbb{E} [\pi(0,z)\prod_{z'\in U}\pi
(0,z')^{n_{o,z'}} ]}{\mathbb{E} [\prod_{z'\in U}\pi
(0,z')^{n_{o,z'}} ]},
\end{equation}
by Bayes' rule. 

Consider the map $T^*\dvtx W_\infty^{\mathrm{tr}}\to W_\infty^{\mathrm{tr}}$
that takes $ (x_i )_{i\leq0}$ to $ (x_i-x_{-1}
)_{i\leq-1}$. Let $\mathcal{I}$ be the set of probability measures on
$W_\infty^{\mathrm{tr}}$ that are invariant under $T^*$ and $\mathcal
{E}$ be the set of extremal points of $\mathcal{I}$. Each $\alpha\in
\mathcal{I}$ (resp., $\alpha\in\mathcal{E})$ corresponds to a transient
process with stationary (resp., stationary and ergodic) increments and
induces a probability measure $Q_{\alpha}$ on particle paths $
(X_i )_{i\in\mathbb{Z}}$. The associated \textit{mean drift} is
$m(\alpha) := \int(x_o-x_{-1} )\,\dd\alpha=Q_{\alpha
}(X_1-X_o)$. Define
%
%
\begin{equation}\label{parisolur}\quad
Q_{\alpha}^w(\cdot):=Q_{\alpha}\bigl( \cdot |\sigma(X_i\dvtx i\leq0)\bigr)(w)
\quad\mbox{and}\quad q_{\alpha}(w,z):=Q_{\alpha}^w(X_1=z)
\end{equation}
for $\alpha$-a.e. $w$ and $z\in U$. 

With this notation,
%
%
\begin{equation}\label{divaneasik}
I_a(\xi)=\mathop{\inf_{\alpha\in\mathcal{E}\dvt}}_{m(\alpha)=\xi
}\mathfrak
{I}_a(\alpha)
\end{equation}
for every $\xi\neq0$, where
%
%
\begin{equation}\label{hayirsh}
\mathfrak{I}_a(\alpha):=\int_{W_\infty^{\mathrm{tr}}} \biggl[\sum_{z\in
U}q_\alpha(w,z)\log\frac{q_\alpha(w,z)}{q(w,z)} \biggr] \,\dd\alpha(w).
\end{equation}

Rassoul-Agha \cite{FirasLDP04} generalizes Varadhan's result to a class
of mixing environments and also to some other models of random walk on
$\mathbb{Z}^d$.

In Section \ref{missinglink}, we will summarize the known qualitative
properties of $I_a$. In particular, we will state some regularity
results which are valid under a certain transience condition of
Sznitman. The next subsection is devoted to introducing this condition,
which involves what are called \textit{regeneration times}.
%

\subsection{Regeneration times and Sznitman's condition}

Take a unit vector $\hat{u}\in\mathcal{S}^{d-1}$.
Define a sequence $ (\tau_m )_{m\geq0}= (\tau_m(\hat
{u}) )_{m\geq0}$ of random times, which are referred to as
\textit{regeneration times} (relative to $\hat{u}$), by $\tau_o:=0$ and
%
%
\begin{eqnarray}\label{regenerationtimes_m}
\tau_{m}&:=&\inf\{j>\tau_{m-1}\dvtx\langle X_i,\hat{u}\rangle<\langle
X_j,\hat{u}\rangle\leq\langle X_k,\hat{u}\rangle\nonumber\\[-8pt]\\[-8pt]
&&\hspace*{76pt}\mbox{for all }i,k\mbox
{ with }i<j<k \}\nonumber
\end{eqnarray}
for every $m\geq1$.
(Regeneration times first appeared in the work of Kesten
\cite{KestenReg} on one-dimensional RWRE. They were adapted to the
multidimensional setting
by Sznitman and Zerner; see \cite{SznitmanZerner99}.)
If the walk is directionally transient relative to $\hat{u}$, that is, if
%
%
\begin{equation}\label{transience}
P_o \Bigl(\lim_{n\to\infty}\langle X_n,\hat{u}\rangle=\infty\Bigr)=1,
\end{equation}
then $P_o (\tau_m<\infty)=1$ for every $m\geq1$.
As shown in \cite{SznitmanZerner99}, the significance of $ (\tau
_m )_{m\geq1}$ is due to the fact that
\[
(X_{\tau_m+1}-X_{\tau_m},X_{\tau_m+2}-X_{\tau_m},\ldots,X_{\tau
_{m+1}}-X_{\tau_m}, \tau_{m+1}-\tau_m )_{m\geq1}
\]
is an i.i.d. sequence under $P_o$ when $\omega=(\omega_x)_{x\in\mathbb
{Z}^d}$ is an i.i.d. collection.

The walk is said to satisfy Sznitman's transience condition ($\mathbf
{T},\hat{u}$)
if (\ref{transience}) holds and
%
%
\begin{equation}\label{moment}
E_o \Bigl[\sup_{1\leq i\leq\tau_1(\hat{u})}\exp\{c |X_i
| \} \Bigr]<\infty\qquad\mbox{for some }c>0.
\end{equation}

Define the \textit{first backtracking time} of the walk to be
%
%
\begin{equation}\label{sisede}
\beta=\beta(\hat{u}):=\inf\{i\geq0\dvtx\langle X_i,\hat{u}\rangle
<\langle X_o,\hat{u}\rangle\}.
\end{equation}
The following lemmas list some important facts regarding regenerations.
\begin{lemma}[(Sznitman \cite{SznitmanT})]\label{madabir}
Assume $d\geq2$, (\ref{ellipticity}), (\ref{i.i.d.}) and that
($\mathbf{T},\hat{u}$) holds for some $\hat{u}\in\mathcal{S}^{d-1}$.
Then:
\begin{enumerate}[(a)]
\item[(a)] $P_o(\beta(\hat{u})=\infty)>0$ and $\tau_1(\hat{u})$ has
finite $P_o$-moments of arbitrary order;
\item[(b)] the LLN holds with a limiting velocity $\xi_o$ such that
$\langle\xi_o,\hat{u}\rangle>0$;
\item[(c)] ($\mathbf{T},\hat{v}$) is satisfied for every $\hat{v}\in
\mathcal{S}^{d-1}$ such that $\langle\xi_o,\hat{v}\rangle>0$.
\end{enumerate}
\end{lemma}
\begin{lemma}\label{madaiki}
Assume (\ref{ellipticity}) and (\ref{i.i.d.}). If the walk is
nonnestling (see Definition~\ref{nonnennon}) relative to some $\hat
{u}\in\mathcal{S}^{d-1}$, then
%
%
\begin{equation}\label{cisembalkan}
E_o [\exp\{c\tau_1(\hat{u}) \} ]<\infty
\end{equation}
for some $c>0$. In particular, ($\mathbf{T},\hat{u}$) is satisfied. On
the other hand, if the walk is nestling, then (\ref{cisembalkan}) fails
to hold for every $\hat{u}\in\mathcal{S}^{d-1}$ and $c>0$.
\end{lemma}
\begin{pf}
The first statement is proved in \cite{SznitmanSlowdown}. The second
statement follows immediately from the fact that $I_a(0)=0$ when the
walk is nestling (see \cite{Varadhan03}).
\end{pf}
\begin{lemma}\label{madauc}
Assume (\ref{ellipticity}) and (\ref{i.i.d.}). If the walk is
nonnestling and some $\hat{v}\in\mathcal{S}^{d-1}$ satisfies $\langle
\xi_o,\hat{v}\rangle>0$, then
\[
E_o [\exp\{c\tau_1(\hat{v}) \} ]<\infty
\]
for some $c>0$.
\end{lemma}
\begin{pf}
This is Lemma 8 of \cite{YilmazQequalsA}.
\end{pf}
\begin{corollary}\label{mazaltov}
Assume $d\geq2$, (\ref{ellipticity}), (\ref{i.i.d.}) and that ($\mathbf
{T},\hat{u}$)
holds for some $\hat{u}\in\mathcal{S}^{d-1}$. Since $\xi
_o\neq0$, there exists a $z\in U$ such that $\langle\xi_o,z\rangle>0$. Then:
\begin{enumerate}[(a)]
\item[(a)] $P_o(\beta(z)=\infty)>0$ and $\tau_1(z)$ has finite
$P_o$-moments of arbitrary order;
\item[(b)] if the walk is nonnestling, then there exists a $c_1>0$
such that
\[
E_o [\exp\{2c_1\tau_1(z) \} ]<\infty;
\]
\item[(c)] if the walk is nestling, then there exists a $c_1>0$ such
that
\[
E_o \Bigl[\sup_{1\leq i\leq\tau_1(z)}\exp\{c_1 |X_i |
\} \Bigr]<\infty.
\]
\end{enumerate}
\end{corollary}
%

\subsection{Qualitative properties of the quenched and the averaged
rate functions}\label{missinglink}

Denote the zero-sets of $I_q$ and $I_a$ by $\mathcal{N}_q:= \{\xi\in
\mathbb{R}^d\dvtx I_q(\xi)=0 \}$ and $\mathcal{N}_a:= \{\xi\in
\mathbb{R}^d\dvtx I_a(\xi)=0 \}$, respectively. The following theorem
summarizes some of the known qualitative properties of the quenched and
the averaged rate functions when $d\geq2$. The rest of the known
properties are given in Section \ref{cizbord}.
\begin{theorem}\label{previousqual}
Assume $d\geq2$, (\ref{ellipticity}) and (\ref{i.i.d.}).
Then:
\begin{enumerate}[(a)]
\item[(a)] $I_q$ and $I_a$ are convex, $I_q(0)=I_a(0)$ and $\mathcal
{N}_q=\mathcal{N}_a$ (see \cite{Varadhan03});
\item[(b)] if the walk is nonnestling, then:
\begin{enumerate}[(ii)]
\item[(i)] $\mathcal{N}_a$ consists of the true velocity $\xi_o$ (see
\cite{Varadhan03});
\item[(ii)] $I_a$ is strictly convex and analytic on an open set
$\mathcal{A}_a$ containing $\xi_o$ (see \cite{JonOferLDP08,YilmazAveraged});
\end{enumerate}
\item[(c)] if the walk is nestling, then $\mathcal{N}_a$ is a line
segment containing the origin that can extend in one or both directions
(see \cite{Varadhan03}); it cannot extend in both directions when
$d=2$ (see \cite{ZernerMerkl}) or when $d\geq5$ (see \cite{BergerZeroSet});
\item[(d)] if the walk is nestling, but ($\mathbf{T},\hat{u}$) is
satisfied for some $\hat{u}\in\mathcal{S}^{d-1}$, then:
\begin{enumerate}[(iii)]
\item[(i)]  the origin is an endpoint of $\mathcal{N}_a$ (see \cite{SznitmanT});
\item[(ii)] $I_a$ is strictly convex and analytic on an open set
$\mathcal{A}_a$ (see \cite{YilmazAveraged});
\item[(iii)]  there exists a $(d-1)$-dimensional smooth surface patch
$\mathcal{A}_a^b$ such that $\xi_o\in\mathcal{A}_a^b\subset\partial
\mathcal{A}_a$ (see \cite{YilmazAveraged});
\item[(iv)]  the unit vector $\eta_o$ normal to $\mathcal{A}_a^b$ (and
pointing in $\mathcal{A}_a$) at $\xi_o$ satisfies $\langle\eta_o,\xi
_o\rangle>0$ (see \cite{YilmazAveraged});
\item[(v)] $I_a(t\xi)=tI_a(\xi)$ for every $\xi\in\mathcal{A}_a^b$ and
$t\in[0,1]$ (see \cite{JonOferLDP08}).
\end{enumerate}
\end{enumerate}
\end{theorem}
%

\subsection{Comparing the quenched and the averaged rate
functions}\label{cizbord}

Assume (\ref{ellipticity}) and (\ref{i.i.d.}). It is clear that
\begin{eqnarray*}
\mathcal{D}:\!&=& \{(\xi_1,\ldots,\xi_d)\in\mathbb{R}^d\dvtx|\xi_1|+\cdots
+|\xi_d|\leq1 \}= \{\xi\in\mathbb{R}^d\dvtx I_a(\xi)<\infty\}
\\
&=& \{\xi\in\mathbb{R}^d\dvtx I_q(\xi)\leq-\log\delta\}.
\end{eqnarray*}
For any $\xi\in\mathbb{R}^d$, $I_a(\xi)\leq I_q(\xi)$ by Jensen's
inequality and Fatou's lemma. Moreover, when the support of $\mathbb
{P}$ is not a singleton, $I_a<I_q$ at some interior points of $\mathcal
{D}$ (see Proposition 4 of \cite{YilmazQequalsA}).

The following theorem considers ballistic walks in dimensions four and
higher, and says that the quenched and the averaged rate functions are
identically equal on a set whose interior contains $\mathcal
{N}_a\setminus\{0\}$.
\begin{theorem}[(Yilmaz \cite{YilmazQequalsA})]\label{QequalsA}
Assume $d\geq4$, (\ref{ellipticity}), (\ref{i.i.d.}) and that
($\mathbf{T},\hat{u}$) holds for some $\hat{u}\in\mathcal{S}^{d-1}$.
Then:
\begin{enumerate}[(a)]
\item[(a)] if the walk is nonnestling, $I_q=I_a$ on an open set
$\mathcal{A}_{\eq}$ containing $\xi_o$;
\item[(b)] if the walk is nestling:
\begin{enumerate}[(iii)]
\item[(i)] $I_q=I_a$ on an open set $\mathcal{A}_{\eq}$;
\item[(ii)]  there exists a $(d-1)$-dimensional smooth surface patch
$\mathcal{A}_{\eq}^b$ such that $\xi_o\in\mathcal{A}_{\eq}^b\subset
\partial\mathcal{A}_{\eq}$;
\item[(iii)]  the unit vector $\eta_o$ normal to $\mathcal{A}_{\eq}^b$
(and pointing in $\mathcal{A}_{\eq}$) at $\xi_o$ satisfies $\langle\eta
_o,\xi_o\rangle>0$;
\item[(iv)] $I_q(t\xi)=tI_q(\xi)=tI_a(\xi)=I_a(t\xi)$ for every $\xi\in
\mathcal{A}_{\eq}^b$ and $t\in[0,1]$.
\end{enumerate}
\end{enumerate}
\end{theorem}

Assuming $d=1$, (\ref{ellipticity}) and (\ref{i.i.d.}), Comets et al.
\cite{CGZ00} use (\ref{ahubay}) to show that $I_q(\xi)=I_a(\xi)$ if and
only if $\xi=0$ or $I_a(\xi)=0$. In particular, Theorem \ref{QequalsA}
cannot be generalized to $d\geq1$. It turns out that it cannot be
generalized to $d\geq2$ or $3$, either. Indeed, for $d=2,3$, Yilmaz and
Zeitouni \cite{YilmazZeitouni09} provide examples of nonnestling walks
in uniformly elliptic i.i.d. environments for which the quenched and
the averaged rate functions are not identically equal on any open set
containing the true velocity $\xi_o$.
%

\subsection{Dual results for the logarithmic moment generating functions}

For every $\theta\in\mathbb{R}^d$, consider the logarithmic moment
generating functions
%
%
\begin{eqnarray}\label{sisedelal}
\Lambda_q(\theta)&:=&\lim_{n\to\infty}\frac{1}{n}\log E_o^\omega[\exp
\{\langle\theta,X_n\rangle\} ] \quad\mbox{and}\nonumber\\[-8pt]\\[-8pt]
\Lambda_a(\theta)
&:=&\lim_{n\to\infty}\frac{1}{n}\log E_o [\exp\{\langle\theta
,X_n\rangle\} ].\nonumber
\end{eqnarray}
By Varadhan's lemma (see \cite{DemboZeitouniBook}), $\Lambda_q(\theta
)=\sup_{\xi\in\mathbb{R}^d} \{\langle\theta,\xi\rangle- I_q(\xi
) \}=I_q^*(\theta)$, the convex conjugate of $I_q$ at $\theta$.
Similarly, $\Lambda_a(\theta)=I_a^*(\theta)$.

For every $c>0$, define
%
%
\begin{equation}\label{nebiliyim}
\mathcal{C}(c):= \cases{
\{\theta\in\mathbb{R}^d\dvtx|\theta|<c \},&\quad if the walk is
nonnestling,\cr
\{\theta\in\mathbb{R}^d\dvtx|\theta|<c , \Lambda_a(\theta)>0 \}
,&\quad if the walk is nestling.}\hspace*{-28pt}
\end{equation}
In the latter case, $I_a(0)=0$; see Theorem \ref{previousqual}. It
follows from convex duality that
\[
0=I_a(0)=\sup_{\theta\in\mathbb{R}^d} \{\langle\theta,0\rangle-
\Lambda_a(\theta) \}=-\inf_{\theta\in\mathbb{R}^d}\Lambda_a(\theta).
\]
In other words, $\Lambda_a(\theta)\geq0$ for every $\theta\in\mathbb
{R}^d$. The zero-level set $ \{\theta\in\mathbb{R}^d\dvtx\Lambda
_a(\theta)=0 \}$ of the convex function $\Lambda_a$ is convex and
$\mathcal{C}(c)$ is an open ball minus this convex set.

The following theorems state some of the known qualitative properties
of $\Lambda_q$ and $\Lambda_a$.
\begin{theorem}[(Peterson and Zeitouni \cite{JonOferLDP08}, Yilmaz \cite
{YilmazAveraged})]\label{Caveraged}
Assume $d\geq2$, (\ref{ellipticity}) and (\ref{i.i.d.}). Recall (\ref
{nebiliyim}). If ($\mathbf{T},\hat{u}$) holds for some $\hat{u}\in
\mathcal{S}^{d-1}$, then $\Lambda_a$ is analytic on $\mathcal
{C}_a:=\mathcal{C}(c_1)$, where $c_1$ is as in Corollary \ref
{mazaltov}. Moreover, the Hessian $\mathcal{H}_a$ of $\Lambda_a$ is
positive definite on $\mathcal{C}_a$.
\end{theorem}
\begin{theorem}[(Yilmaz \cite{YilmazQequalsA})]\label{Cananbe}
Assume $d\geq4$, (\ref{ellipticity}) and (\ref{i.i.d.}). Recall (\ref
{nebiliyim}). If ($\mathbf{T},\hat{u}$) holds for some $\hat{u}\in
\mathcal{S}^{d-1}$, then there exists a $c_2\in(0,c_1)$ such that
$\Lambda_q=\Lambda_a$ on $\mathcal{C}_{\eq}:=\mathcal{C}(c_2)$.
\end{theorem}

In fact, the regularity properties of $I_a$ that are stated in
Theorem \ref{previousqual} are obtained from Theorem \ref{Caveraged} via convex
duality (see \cite{JonOferLDP08,YilmazAveraged}) and $\mathcal{A}_a=\{
\nabla\Lambda_a(\theta)\dvtx\theta\in\mathcal{C}_a\}$. Similarly, note that
Theorem \ref{QequalsA} is a corollary of Theorem~\ref{Cananbe} and
$\mathcal{A}_{\eq}=\{\nabla\Lambda_a(\theta)\dvtx\theta\in\mathcal
{C}_{\eq}\}$.
%

\section{Our results}\label{inisoyler}

In this paper, we will obtain new results concerning the large
deviation properties of RWRE on $\mathbb{Z}^d$ under the conditions of
Theorems \ref{QequalsA} and~\ref{Cananbe}. In other words, we will
assume that
%
%
\begin{equation}\label{kondisin}
\begin{tabular}{p{323pt}}
\mbox{$d\geq4$, the environment is uniformly elliptic and i.i.d.} \mbox{[see
(\ref{ellipticity}) and (\ref{i.i.d.})]} and ($\mathbf{T},e_1$) holds.
\end{tabular}\hspace*{-32pt}
\end{equation}
Here, we have chosen $e_1$ for convenience. However, there is no loss
of generality, that is, we could have chosen any $\hat{u}\in\mathcal
{S}^{d-1}$; see Lemma \ref{madabir}.
%

\subsection{Existence of harmonic functions: h-transform}

Given any $\theta\in\mathbb{R}^d$, define $\pi^\theta\dvtx\Omega\times
U\to
\mathbb{R}$ by setting
\[
\pi^\theta(\omega,z):=\pi(0,z)\exp\{\langle\theta,z\rangle-\Lambda
_a(\theta)\}
\]
for every $\omega\in\Omega$ and $z\in U$. Our first result concerns the
existence of positive harmonic functions for $\pi^\theta$. (Here, we
use the term \textit{harmonic} in analogy with the continuum case where
$\pi^\theta$ is replaced by a second order elliptic operator.)

\begin{theorem}\label{existence}
Assume (\ref{kondisin}). Recall Theorem \ref{Cananbe}. For every $\theta
\in\mathcal{C}_{\eq}$, there exists an $h(\theta,\cdot)\in L^2(\mathbb
{P})$ such that $\mathbb{P}(h(\theta,\cdot)>0)=1$ and
%
%
\begin{equation}\label{carbonic}\qquad
h(\theta,\omega)=\sum_{z\in U}\pi(0,z)\exp\{\langle\theta,z\rangle
-\Lambda_a(\theta)\}h(\theta,T_z\omega) \qquad\mbox{for $\mathbb
{P}$-a.e. $\omega$.}
\end{equation}
\end{theorem}

Note that $\pi^\theta$ would correspond to a Markov chain on $\Omega$
if $\sum_{z\in U}\pi^\theta(\omega,z)=1$ were true for $\mathbb
{P}$-a.e. $\omega$. However, as we will see, the latter condition is
not satisfied unless $\theta=0$. Nevertheless, (\ref{carbonic}) enables
us to define an environment kernel (as in Definition \ref{ortamkeli})
related to $\pi^\theta$ via the so-called \textit{h-transform}
technique of Doob; see~\cite{PinskyBook}.
\begin{definition}\label{doobabi}
Assume (\ref{kondisin}). For every $\theta\in\mathcal{C}_{\eq}$, define
a new environment kernel $\hat{\pi}^\theta\dvtx\Omega\times U\to\mathbb
{R}^+$ by setting
%
%
\begin{equation}\label{azdinyn}
\hat{\pi}^\theta(\omega,z):=\pi(0,z)\exp\{\langle\theta,z\rangle-\Lambda
_a(\theta)\}\frac{h(\theta,T_z\omega)}{h(\theta,\omega)}
\end{equation}
for every $\omega\in\Omega$ and $z\in U$. This technique is called
\textit{h-transform}.
\end{definition}
%

\subsection{The unique minimizer of Varadhan's variational
formula}\label{unminvar}

Recall the sets $\mathcal{A}_a$ and $\mathcal{A}_{\eq}$ which were
introduced in Theorems \ref{previousqual} and \ref{QequalsA},
respectively. Whenever $\xi\in\mathcal{A}_a$, it is shown in
\cite{YilmazAveraged} that there is a unique minimizer of Varadhan's
variational formula (\ref{divaneasik}) for $I_a(\xi)$. Our second
result reveals the hidden Markovian structure of this minimizer when
(\ref{kondisin}) holds and $\xi\in\mathcal{A}_{\eq}$.

Before stating this theorem, we need to introduce a family of
sub-$\sigma$-algebras of $\mathcal{B}$: for any $\hat{v}\in\mathcal
{S}^{d-1}$ and $n\geq0$, let
%
%
\begin{equation}\label{subsigsub}
\mathcal{B}_n^+(\hat{v}):=\sigma(\omega_x\dvtx\langle x,\hat{v}\rangle
\geq-n).
\end{equation}

\begin{theorem}\label{varminvar}
Assume (\ref{kondisin}). 
Recall Theorem \ref{QequalsA} and Definition \ref{doobabi}. For every
$\xi\in\mathcal{A}_{\eq}$, there exists a unique $\theta\in\mathcal
{C}_{\eq}$ such that $\xi=\nabla\Lambda_a(\theta)$. [By convex duality,
$\theta=\nabla I_a(\xi)$.]
\begin{enumerate}[(a)]
\item[(a)] There exists a unique $\mathbb{Q}_\xi\in M_1(\Omega,\mathcal
{B})$ that satisfies the following:
\begin{enumerate}[(ii)]
\item[(i)] $\mathbb{Q}_\xi$ is $\hat{\pi}^\theta$-invariant, that is,
$\sum_{z\in U}\dd\mathbb{Q}_\xi(T_{-z}\omega)\hat{\pi}^\theta
(T_{-z}\omega,z)=\dd\mathbb{Q}_\xi(\omega)$;
\item[(ii)] $\mathbb{Q}_\xi\ll\mathbb{P}$ on $\mathcal{B}_n^+(e_1)$ for
every $n\geq0$; see (\ref{subsigsub}).
\end{enumerate}
The pair $(\hat{\pi}^\theta,\mathbb{Q}_\xi)$ corresponds to a
stationary Markov chain (with values in $\Omega$) which can be
identified with a $\hat{\mu}_\xi^\infty\in M_1(\Omega\times U^\mathbb
{N})$. The marginal on $\Omega$ of $\hat{\mu}_\xi^\infty$ is $\mathbb
{Q}_\xi$ and $\hat{\pi}^\theta$ is the conditional of $z_1$ given
$\omega$.
\item[(b)] $\hat{\mu}_\xi^\infty$ induces a $\mathbb{Z}^d$-valued
transient process with stationary increments in $U$ via the map
\[
(\omega, z_1,z_2,z_3,\ldots)\mapsto(z_1,z_1+z_2,z_1+z_2+z_3,\ldots).
\]
Extend this process to a probability measure on doubly infinite paths
$(x_i)_{i\in\mathbb{Z}}$ and refer to its restriction to $W_{\infty
}^{\mathrm{tr}}$ as $\mu_\xi^\infty$. With this notation, $\mu_\xi^\infty$ is
the unique minimizer of Varadhan's variational formula (\ref{divaneasik}).
\end{enumerate}
\end{theorem}

In words, when a particle under $P_o$ is conditioned to have asymptotic
mean velocity equal to any given $\xi\in\mathcal{A}_{\eq}$, the
environment Markov chain chooses to switch from its original kernel
$\overline{\pi}$ to the tilted kernel $\hat{\pi}^\theta$ given in (\ref
{azdinyn}), where $\theta=\nabla I_a(\xi)\in\mathcal{C}_{\eq}$. The most
economical tilt in terms of averaged large deviations is realized by an
h-transform.
\begin{remark}
There is an alternative characterization of $\hat{\mu}_\xi^\infty$
[see (\ref{muyucananan})] which involves regeneration times. That
formula (or, rather, its analog for the marginal on $U^\mathbb{N}$ of
$\hat{\mu}_\xi^\infty$) has already appeared in Definition 9 of
\cite{YilmazAveraged}. If one takes (\ref{muyucananan}) as the definition of
$\hat{\mu}_\xi^\infty$, then part (b) of Theorem \ref{varminvar}
becomes essentially a restatement of Theorem 10 of
\cite{YilmazAveraged} (see Theorem \ref{cizme} of the current paper for
details). In other words, the novelty of Theorem \ref{varminvar} lies
in part (a).
\end{remark}
%

\subsection{Equality of the quenched and the averaged minimizers}

The quenched level-3 LDP stated in Theorem \ref{level3LDP} implies the
quenched (level-1) LDP (i.e., Theorem \ref{qLDPgeneric}) via the
contraction principle. Indeed, for any $\xi\in\mathbb{R}^d$, define
%
%
\begin{equation}\label{Axizz}
A_\xi^\infty:= \biggl\{\hat{\alpha}\in M_1(\Omega\times U^\mathbb{N})\dvtx\int
\sum_{(z_i)_{i\geq1}\in U^\mathbb{N}} \dd\hat{\alpha
}(\omega,(z_i)_{i\geq1})z_1=\xi\biggr\}.
\end{equation}
With this notation,
%
%
\begin{equation}\label{level1ratezz}
I_q(\xi)=\inf_{\hat{\alpha}\in A_\xi^\infty}I_{q,3}(\hat{\alpha}).
\end{equation}

Our third result is
as follows.
\begin{theorem}\label{rosminros}
Assume (\ref{kondisin}). For every $\xi\in\mathcal{A}_{\eq}$, the
measure $\hat{\mu}_\xi^\infty$ (which is obtained in Theorem
\ref{varminvar}) is the unique minimizer of (\ref{level1ratezz}).
\end{theorem}

We already know from Theorem \ref{QequalsA} that the quenched and the
averaged rate functions $I_q$ and $I_a$ are equal on $\mathcal
{A}_{\eq}$. The natural interpretation of Theorem \ref{rosminros} is
that, for $\mathbb{P}$-a.e. $\omega$, when a particle under $P_o^\omega
$ is conditioned to have asymptotic mean velocity equal to any given
$\xi\in\mathcal{A}_{\eq}$, the environment Markov chain chooses to
switch from its original kernel $\overline{\pi}$ to the tilted kernel
$\hat{\pi}^\theta$. Compare this with the last paragraph of the
previous subsection.

Since the contraction from level-3 to level-1 may be done in two steps
(instead of one), the following is an immediate consequence of Theorem
\ref{rosminros}.

\begin{corollary}\label{duzbac}
Assume (\ref{kondisin}). For every $\xi\in\mathcal{A}_{\eq}$, let $\hat
{\mu}_\xi\in M_1(\Omega\times U)$ be the marginal
of $\hat{\mu}_\xi^\infty\in M_1(\Omega\times U^\mathbb{N})$. With this
notation, $\hat{\mu}_\xi$ is the unique minimizer of the variational
formula
\[
I_q(\xi)=\inf_{\mu\in A_\xi} \mathfrak{I}_q^{**}(\mu)
\]
given in (\ref{level1rate}).
\end{corollary}
%

\subsection{Modifying Rosenbluth's variational formula}

Recall Rosenbluth's variational formula
\[
I_q(\xi)=\inf_{\mu\in A_\xi} \mathfrak{I}_q(\mu)
\]
given in (\ref{level1ratetilde}). Its advantage over (\ref{level1rate})
is that $\mathfrak{I}_q$ has a simple formula, whereas $\mathfrak
{I}_q^{**}$ does not.
Corollary \ref{duzbac} identifies the unique minimizer of (\ref
{level1rate}) when (\ref{kondisin}) holds and $\xi\in\mathcal{A}_{\eq}$.
We would like to obtain an analogous result for (\ref
{level1ratetilde}). However, as we illustrate below, there is a problem.

We express Rosenbluth's formula in the following way:
%
%
\begin{equation}\label{guzelyaziss}
I_q(\xi)=\inf\{H(\mu)\dvtx\mu\in A_\xi\cap M_1'(\Omega\times U)\},
\end{equation}
where
%
%
\begin{equation}\label{halitakca}
H(\mu):=\int_{\Omega}\sum_{z\in U} \dd\mu(\omega,z)\log\frac
{\dd\mu(\omega,z)}{\dd(\mu)^1(\omega)\pi(0,z)}
\end{equation}
denotes relative entropy and
\begin{eqnarray*}
M_1'(\Omega\times U)&:=&\biggl\{\mu\in M_1(\Omega\times U)\dvtx(\mu)^1=(\mu
)^2\ll\mathbb{P},\\
&&\hspace*{34.6pt} \frac{\dd\mu(\cdot,z)}{\dd(\mu)^1(\cdot
)}>0 \mbox{ for every }z\in U \biggr\}.
\end{eqnarray*}
In light of Corollary \ref{duzbac}, a natural minimizer candidate
for (\ref{guzelyaziss}) is $\hat{\mu}_\xi$. Note that $\hat{\mu}_\xi$
is an element of $M_1'(\Omega\times U)$ if and only if its marginal
$\mathbb{Q}_\xi$ is absolutely continuous relative to $\mathbb{P}$ on
$\mathcal{B}$. However, all we know is that $\mathbb{Q}_{\xi}\ll\mathbb
{P}$ on $\mathcal{B}_n^+(e_1)$ for every $n\geq0$; see Theorem \ref{varminvar}.

Instead of trying to show that $\hat{\mu}_\xi$ is an element of
$M_1'(\Omega\times U)$, we will replace $M_1'(\Omega\times U)$ by a
larger set that contains $\hat{\mu}_\xi$.
\begin{definition}\label{omurtukendi}
A measure $\mu\in M_1(\Omega\times U)$ is said to be in $M_1''(\Omega
\times U)$ if it satisfies the following conditions:
\begin{enumerate}[(a)]
\item[(a)] $(\mu)^1=(\mu)^2$; see (\ref{albaytey});
\item[(b)] $\hat{\pi}(\cdot,z):=\frac{\dd\mu(\cdot,z)}{\dd(\mu)^1(\cdot
)}>0$ for every $z\in U$;
\item[(c)] there exists a $\hat{v}\in\mathcal{S}^{d-1}$ such that
$P_o^{\hat{\pi}}(\lim_{n\to\infty}\langle X_n,\hat{v}\rangle=\infty
)=1$; see Definition \ref{ortamkeli};
\item[(d)] $(\mu)^1\ll\mathbb{P}$ on $\mathcal{B}_n^+(\hat{v})$ for
every $n\geq0$; see (\ref{subsigsub}).
\end{enumerate}
\end{definition}
\begin{theorem}\label{gobenk}
Assume (\ref{ellipticity}). Recall (\ref{halitakca}) and Definition \ref
{omurtukendi}. For every \mbox{$\xi\neq0$},
%
%
\begin{equation}\label{guzelyazisss}
I_q(\xi)=\inf\{H(\mu)\dvtx\mu\in A_\xi\cap M_1''(\Omega\times U)\}.
\end{equation}
\end{theorem}

Our last result is
the following theorem.
\begin{theorem}\label{tezbitti}
Assume (\ref{kondisin}). For every $\xi\in\mathcal{A}_{\eq}$, $\hat{\mu
}_\xi$ is the unique minimizer of (\ref{guzelyazisss}).
\end{theorem}

Note that (\ref{guzelyazisss}) does not involve any complex conjugation
and, therefore, is simpler (i.e., more explicit) than (\ref
{level1rate}). Because of this, we believe that Theorem \ref{tezbitti}
is more useful than Corollary \ref{duzbac}.

%

\subsection{Some questions and comments}

\begin{enumerate}
\item When (\ref{kondisin}) holds and $\xi\in\mathcal{A}_{\eq}$,
Theorem \ref{varminvar} states that $\mathbb{Q}_{\xi}\ll\mathbb{P}$ on
$\mathcal{B}_n^+(e_1)$ for every $n\geq0$. On the other hand, it is not
known if $\mathbb{Q}_{\xi}\ll\mathbb{P}$ on $\mathcal{B}$. Is the
latter statement true? Note that, when $\xi=\xi_o$, this question is of
great interest (in its own right) because $\mathbb{Q}_{\xi_o}$ is the
invariant measure from the point of view of the particle. 

Bolthausen and Sznitman \cite{BolthausenSznitman02} prove that $\mathbb
{Q}_{\xi}\ll\mathbb{P}$ on $\mathcal{B}$ when $\xi=\xi_o$ and the
disorder in the environment is low. One expects their argument to work
when $|\xi-\xi_o|$ is small. However, their technique does not
generalize to the case where the disorder is not low.
\item The limitation of our results is that they are valid when
(\ref{kondisin}) holds and $\xi\in\mathcal{A}_{\eq}$, and their proofs
break down if any of these assumptions are weakened. Therefore, it is
natural to ask the following question: in the context of
multidimensional RWRE, does the connection between h-transform and
large deviations exist under more general conditions? 
Note that such a connection has been established (i) for walks with
bounded jumps on $\mathbb{Z}$ in stationary and ergodic environments
(see \cite{YilmazQuenched}), and (ii) for space--time walks in
dimensions $3+1$ and higher (see \cite{YilmazSpaceTime}).
\end{enumerate}

The rest of this paper is devoted to the proofs of our results. Most of
our efforts are focused on Theorems \ref{existence} and \ref
{varminvar}, which are established in Sections \ref{kangibir} and \ref
{kangiiki}, respectively. The remaining results (i.e., Theorems \ref
{rosminros}, \ref{gobenk} and \ref{tezbitti}) are obtained in Section
\ref{kangiuc}.
%

\section{Proof of the existence of harmonic functions}\label{kangibir}

\subsection{An $L^2$ estimate}

Assume (\ref{kondisin}). Recall (\ref{sisede}) and (\ref{sisedelal}).
For every $n\geq1$, $\theta\in\mathbb{R}^d$ and $\omega\in\Omega$, define
%
%
\begin{eqnarray}
\label{kafes}
g_n(\theta,\omega)&:=&E_o^\omega[\exp\{\langle\theta
,X_{H_n}\rangle-\Lambda_a(\theta)H_n \},H_n=\tau_k\nonumber\\[-8pt]\\[-8pt]
&&\hspace*{74.8pt} \mbox{for some
}k\geq1, \beta=\infty]\nonumber
\end{eqnarray}
and
\begin{equation}
\label{kameriye}\quad
h_n(\theta,\omega):=E_o^\omega[\exp\{\langle\theta
,X_{H_n}\rangle-\Lambda_a(\theta)H_n \},H_n=\tau_k \mbox{ for some
}k\geq1 ],
\end{equation}
where
\[
H_n:=\inf\{i\geq0\dvtx\langle X_i,e_1\rangle\geq n \}.
\]
\begin{lemma}[(Yilmaz \cite{YilmazQequalsA})]\label{eskiL2lemma}
Assume (\ref{kondisin}). There exists a $c_2\in(0,c_1)$ such that
%
%
\begin{equation}
\label{lemsambir}
\liminf_{n\to\infty}\mathbb{E} \{g_n(\theta,\cdot) \}>0
\end{equation}
and
\begin{equation}
\label{lemsamiki}
\sup_{n\geq1}\mathbb{E} \{g_n(\theta,\cdot)^2 \}<\infty
\end{equation}
for every $\theta\in\mathcal{C}_{\eq}:=\mathcal{C}(c_2)$; see (\ref{nebiliyim}).
\end{lemma}
\begin{pf}
This constitutes the core of the proof of Theorem \ref{Cananbe}. For
the convenience of the reader, we will give a sketch of the argument.
See Lemmas 11 and~12 of \cite{YilmazQequalsA} for the complete proof.

It is shown in Lemma 12 of \cite{YilmazAveraged} that for every $\theta
\in\mathcal{C}_a:=\mathcal{C}(c_1)$,
%
%
\begin{equation}\label{yalebbim}
E_o [\exp\{\langle\theta,X_{\tau_1}\rangle-\Lambda_a(\theta
)\tau_1\} |\beta=\infty]=1.
\end{equation}
[Note that Theorem \ref{Caveraged} follows from (\ref{yalebbim}) by the
implicit function theorem.]
For every $y\in\mathbb{Z}^d$, let
\[
q^\theta(y):=E_o [\exp\{\langle\theta,X_{\tau_1}\rangle
-\Lambda_a(\theta)\tau_1\}, X_{\tau_1}=y |\beta=\infty].
\]
Since $\sum_{y\in\mathbb{Z}^d}q^\theta(y)=1$ by (\ref{yalebbim}), $
(q^\theta(y) )_{y\in\mathbb{Z}^d}$ defines a random walk
$(Y_k)_{k\geq0}$ on~$\mathbb{Z}^d$. For every $n\geq1$, $\mathbb{E}
\{g_n(\theta,\cdot) \}/P_o(\beta=\infty)$ is equal to the
probability of the event $\{\langle Y_k,e_1\rangle=n\mbox{ for some
}k\geq1\}$. By renewal theory, this probability is easily shown to
converge to a nonzero limit. In particular, (\ref{lemsambir}) follows.

For every $x,\tilde{x}\in\mathbb{Z}^d$, consider two independent walks
$X=X(x):=(X_i)_{i\geq0}$ and $\tilde{X}=\tilde{X}(\tilde{x}):=(\tilde
{X}_j)_{j\geq0}$, starting at $x$ and $\tilde{x}$, respectively, in the
same environment. Denote their joint quenched law and joint averaged
law by $P_{x,\tilde{x}}^\omega:=P_x^\omega\otimes P_{\tilde
{x}}^\omega$ and $P_{x,\tilde{x}}(\cdot):=\mathbb{E}\{P_{x,\tilde
{x}}^\omega(\cdot)\}$, respectively. As usual, $E_{x,\tilde{x}}^\omega$
and $E_{x,\tilde{x}}$ refer to expectations under $P_{x,\tilde
{x}}^\omega$ and $P_{x,\tilde{x}}$, respectively.

Clearly, $P_{x,\tilde{x}}\neq P_x \otimes P_{\tilde{x}}$. On the
other hand, the two walks do not know that they are in the same
environment unless their paths intersect. In particular, for any event
$A$ involving $X$ and $\tilde{X}$,
%
%
\begin{equation}\label{cokkibar}
P_{x,\tilde{x}}(A\cap\{\gamma_1=\infty\})=P_x \otimes P_{\tilde
{x}}(A\cap\{\gamma_1=\infty\}),
\end{equation}
if $x\neq\tilde{x}$, where
%
%
\begin{equation}\label{hatirlatmak}\quad
\gamma_1:=\inf\{m\in\mathbb{Z}\dvtx X_i=\tilde{X}_j\mbox{ for some
}i\geq
0, j\geq0\mbox{, and }\langle X_i,e_1\rangle=m\}.
\end{equation}

Similar to the random times 
$(H_n)_{n\geq0}$ and $\beta$ for $X$, define 
$(\tilde{H}_n)_{n\geq0}$ and $\tilde{\beta}$ for $\tilde{X}$.
The proof of (\ref{lemsamiki}) makes use of the \textit{joint
regeneration levels} of $X$ and $\tilde{X}$, which are elements of
\[
\mathcal{L}:=\{n\geq0\dvtx\langle X_i,e_1\rangle\geq n\mbox{ and
}\langle
\tilde{X}_j,e_1\rangle\geq n\mbox{ for every }i\geq H_n\mbox{ and }j\geq
\tilde{H}_n\}.
\]
%
Note that if the starting points $x$ and $\tilde{x}$ are both in
$\mathbb{V}_d:= \{z\in\mathbb{Z}^d\dvtx\langle z,e_1\rangle=0 \}
$, then
\[
0\in\mathcal{L} \quad\iff\quad \beta=\tilde{\beta}=\infty\quad\iff\quad
l_1:=\inf\mathcal{L}=0.
\]

For every $n\geq1$ and $\theta\in\mathcal{C}_a$, define
\[
f(\theta,n,X,\tilde{X}):=\exp\{\langle\theta,X_{H_n}\rangle-\Lambda
_a(\theta)H_n \}\exp\{\langle\theta,\tilde{X}_{\tilde{H}_n}\rangle
-\Lambda_a(\theta)\tilde{H}_n\}.
\]
With this notation,
%
%
\begin{equation}\label{gizyuc}
\mathbb{E} \{g_n(\theta,\cdot)^2 \}=E_{o,o} [f(\theta
,n,X,\tilde{X}), n\in\mathcal{L}, l_1=0 ].
\end{equation}

By Lemma \ref{nonintersectlemma} (stated below), the random paths $X$
and $\tilde{X}$ intersect finitely many times and the probability that
they intersect far away from the origin is exponentially small.
Conditioned on the first joint regeneration level after the last
intersection, the right-hand side of (\ref{gizyuc}) can be written as a
product of two terms. The first term is shown to be finite, by renewal
theory, when $\theta\in\mathcal{C}(c_2)$ with a small enough $c_2\in
(0,c_1)$, and the second term is bounded from above by $\mathbb{E}
\{g_n(\theta,\cdot) \}^2\leq1$ since the walks can be thought of
as taking place in independent environments.
\end{pf}

As mentioned in the sketch above, the following lemma is central to the
proof of (\ref{lemsamiki}).
\begin{lemma}\label{nonintersectlemma}
Assume (\ref{kondisin}). Recall (\ref{hatirlatmak}) and let $\mathbb
{V}_d':=\mathbb{V}_d\setminus\{0\}$.
Then
%
%
\begin{equation}\label{eyluldebasvur}
\inf_{z\in\mathbb{V}_d'}P_{o,z}(l_1=0)\geq\inf_{z\in\mathbb
{V}_d'}P_{o,z}(\gamma_1=\infty, l_1=0)>0.
\end{equation}
\end{lemma}
\begin{pf}
Assume (\ref{kondisin}). We saw in part (a) of Corollary \ref{mazaltov}
that $\tau_1$ has finite moments of arbitrary order. Therefore, the
second inequality follows from Proposition 3.1 (for $d\geq5$) and
Proposition 3.4 (for $d=4$) of the recent work of Berger and Zeitouni
\cite{NoamOferCLT}. (The proofs of these propositions are based on
certain Green's function estimates which fail to hold unless $d\geq4$.)
Since the first inequality is clear, we are done.
\end{pf}
\begin{remark}
It is easy to see that the first infimum in (\ref{eyluldebasvur}) is
positive when $d=2,3$ as well. However, we will not need this fact in
what follows.
\end{remark}
\begin{lemma}\label{yeniL2lemma}
Assume (\ref{kondisin}). Recall (\ref{kameriye}). For every $\theta\in
\mathcal{C}_{\eq}$,
%
%
\begin{equation}
\label{lemsambin}
\liminf_{n\to\infty}\mathbb{E} \{h_n(\theta,\cdot) \}>0
\end{equation}
and
\begin{equation}
\label{lemsamikibin}
\sup_{n\geq1}\mathbb{E} \{h_n(\theta,\cdot)^2 \}<\infty.
\end{equation}
\end{lemma}
\begin{pf}
Recall the notation in the sketch of the proof of Lemma \ref
{eskiL2lemma}. By definition, $h_n(\theta,\omega)\geq g_n(\theta,\omega
)$ for every $n\geq1$, $\theta\in\mathcal{C}_{\eq}$ and $\omega\in\Omega
$. Hence, (\ref{lemsambin}) follows immediately from (\ref{lemsambir}).

For every $n\geq1$ and $\theta\in\mathcal{C}_{\eq}$,
%
%
\begin{eqnarray}\qquad
\mathbb{E} \{h_n(\theta,\cdot)^2 \}&=&E_{o,o} [f(\theta
,n,X,\tilde{X}), n\in\mathcal{L} ]\nonumber\\
&=&\sum_{k=0}^{n}\sum_{z\in\mathbb{V}_d}E_{o,o} [f(\theta,n,X,\tilde
{X}),l_1=k,\tilde{X}_{\tilde{H}_k}-X_{H_k}=z,n\in\mathcal{L}
]\nonumber\\
\label{cokheyecanli}
&=&\sum_{k=0}^{n}\sum_{z\in\mathbb{V}_d}E_{o,o} [f(\theta,k,X,\tilde
{X}),l_1=k,\tilde{X}_{\tilde{H}_k}-X_{H_k}=z ]\\
&&\hspace*{34.5pt}{}\times\ee^{-\langle\theta,z\rangle}E_{o,z}
[f(\theta,n-k,X,\tilde{X}),n-k\in\mathcal{L} |l_1=0
]\nonumber\\
\label{olacakbence}
&\leq& E_{o,o} [f(\theta,l_1,X,\tilde{X}) ] \Bigl(\inf_{z\in
\mathbb{V}_d}P_{o,z}(l_1=0) \Bigr)^{-1}\nonumber\\[-8pt]\\[-8pt]
&&{}\times\mathop{\sup_{0\leq k\leq n}}_{z\in
\mathbb{V}_d}\mathbb{E} \{g_{n-k}(\theta,\cdot)g_{n-k}(\theta
,T_z\cdot) \}\nonumber\\
\label{madlenciko}
&\leq& E_{o,o} [f(\theta,l_1,X,\tilde{X}) ] \Bigl(\inf_{z\in
\mathbb{V}_d}P_{o,z}(l_1=0) \Bigr)^{-1}\nonumber\\[-8pt]\\[-8pt]
&&{}\times\sup_{m\geq1}\mathbb{E} \{
g_m(\theta,\cdot)^2 \}.\nonumber
\end{eqnarray}
Indeed, we have (\ref{cokheyecanli}) by the independence structure
which is still valid for common regeneration blocks. (\ref
{olacakbence}) follows by noting that
\begin{eqnarray*}
&&\ee^{-\langle\theta,z\rangle}E_{o,z} [f(\theta
,n-k,X,\tilde{X}),n-k\in\mathcal{L} |l_1=0 ]\\
&&\qquad=\bigl(P_{o,z}(l_1=0) \bigr)^{-1}\mathbb{E} \{g_{n-k}(\theta,\cdot
)g_{n-k}(\theta,T_z\cdot) \}\\
&&\qquad\leq\Bigl(\inf_{z\in\mathbb{V}_d}P_{o,z}(l_1=0) \Bigr)^{-1}
\mathop{\sup_{0\leq k\leq n}}_{z\in\mathbb{V}_d}\mathbb{E} \{
g_{n-k}(\theta
,\cdot)g_{n-k}(\theta,T_z\cdot) \}.
\end{eqnarray*}
The third term in (\ref{madlenciko}) is obtained using the Schwarz
inequality and it is finite by Lemma \ref{eskiL2lemma}. Similarly, the
second term in (\ref{madlenciko}) is finite by Lemma \ref
{nonintersectlemma}. Therefore, to prove (\ref{lemsamikibin}), it
suffices to show that the first term in (\ref{madlenciko}) is also finite.

By H\"older's inequality,
%
%
\begin{eqnarray}\label{negobe}\quad
&&E_{o,o} [f(\theta,l_1,X,\tilde{X}) ]\nonumber\\
&&\qquad=\sum_{k=0}^\infty
E_{o,o} [\exp\{\langle\theta,X_{H_k}\rangle-\Lambda_a(\theta
)H_k \}\nonumber\\
&&\qquad\quad\hspace*{36.3pt}{}\times\exp\{\langle\theta,\tilde{X}_{\tilde{H}_k}\rangle-\Lambda
_a(\theta)\tilde{H}_k\}, l_1=k ]\nonumber\\
&&\qquad \leq\sum_{k=0}^\infty E_{o,o} [\exp\{4\langle\theta
,X_{H_k}\rangle-4\Lambda_a(\theta)H_k \} ]^{1/4}\nonumber\\
&&\qquad\quad\hspace*{13.7pt}{}\times
E_{o,o}[\exp\{4\langle\theta,\tilde{X}_{\tilde{H}_k}\rangle-4\Lambda_a(\theta
)\tilde{H}_k\} ]^{1/4}P_{o,o}(l_1=k)^{1/2}\nonumber\\
&&\qquad=\sum_{k=0}^\infty E_o [\exp\{4\langle\theta
,X_{H_k}\rangle-4\Lambda_a(\theta)H_k \}
]^{1/2}P_{o,o}(l_1=k)^{1/2}.
\end{eqnarray}
For any $k\geq1$,
%
%
\begin{eqnarray}
&& E_o [\exp\{4\langle\theta,X_{H_k}\rangle-4\Lambda_a(\theta
)H_k \} ]\nonumber\\
&&\qquad =E_o [\exp\{4\langle\theta,X_{H_k}\rangle-4\Lambda
_a(\theta)H_k \}, H_k\leq\tau_1 ]\nonumber\\
&&\qquad\quad{} +E_o [\exp\{4\langle
\theta,X_{H_k}\rangle-4\Lambda_a(\theta)H_k \},\tau_1<H_k
]\nonumber\\
&&\qquad =E_o [\exp\{4\langle\theta,X_{H_k}\rangle-4\Lambda
_a(\theta)H_k \}, H_k\leq\tau_1 ]\nonumber\\
&&\qquad\quad{}+\sum_{j=1}^{k-1}E_o
[\exp\{4\langle\theta,X_{H_k}\rangle-4\Lambda_a(\theta)H_k \}
,\tau_1=H_j ]\nonumber\\
&&\qquad =E_o [\exp\{4\langle\theta,X_{H_k}\rangle-4\Lambda
_a(\theta)H_k \}, H_k\leq\tau_1 ]\nonumber\\
&&\qquad\quad{} +\sum_{j=1}^{k-1}E_o [\exp\{4\langle\theta,X_{\tau
_1}\rangle-4\Lambda_a(\theta)\tau_1 \},\tau_1=H_j ]\nonumber\\
&&\qquad\quad\hspace*{26.4pt}{}\times E_o
[\exp\{4\langle\theta,X_{H_{k-j}}\rangle-4\Lambda_a(\theta
)H_{k-j} \} |\beta=\infty]\nonumber\\
&&\qquad \leq E_o \Bigl[\sup_{1\leq n\leq\tau_1}\exp\bigl\{4|\theta
||X_n|-4\bigl(0\wedge\Lambda_a(\theta)\bigr)\tau_1 \bigr\} \Bigr]\nonumber\\
&&\qquad\quad{} \times\Bigl(1+\sup_{1\leq i<k}E_o [\exp\{4\langle
\theta,X_{H_i}\rangle-4\Lambda_a(\theta)H_i \} |\beta=\infty
] \Bigr)\nonumber\\
\label{revirvar}
&&\qquad \leq E_o \Bigl[\sup_{1\leq n\leq\tau_1}\exp\bigl\{4|\theta
||X_n|-4\bigl(0\wedge\Lambda_a(\theta)\bigr)\tau_1 \bigr\} \Bigr] \Bigl(1+\sup
_{1\leq i<k}i\ee^{a_1|\theta|i} \Bigr)\\
\label{cubar}
&&\qquad \leq K_1 \bigl(1+k\ee^{a_1|\theta|k} \bigr).
\end{eqnarray}
Indeed, if the walk is nonnestling, we have $4|\theta||X_n|-4(0\wedge
\Lambda_a(\theta))\tau_1\leq8|\theta|\tau_1$ for every $n\leq\tau_1$.
On the other hand, if the walk is nestling, then $4|\theta
||X_n|-4(0\wedge\Lambda_a(\theta))\tau_1=4|\theta||X_n|$ since $\Lambda
_a(\theta)\geq0$. Therefore, in both cases, the first term in (\ref
{revirvar}) is finite (provided that $4|\theta|<c_1$) and it is denoted
by $K_1$ in (\ref{cubar}). The second term in (\ref{revirvar}) is
obtained using Lemma 28 of \cite{YilmazQequalsA}, where $a_1>0$ is a constant.

It is shown in (the proof of) Lemma 30 of \cite{YilmazQequalsA} that
$E_{o,o} [\ee^{a_3l_1} ]<\infty$ for some $a_3>0$. For
any $k\geq1$,
%
%
\begin{equation}\label{dondolas}
P_{o,o}(l_1=k)\leq E_{o,o}[\ee^{a_3l_1}]\ee^{-a_3k}=:K_2\ee^{-a_3k}.
\end{equation}
Putting (\ref{negobe}), (\ref{cubar}) and (\ref{dondolas}) together, we
conclude that
\begin{eqnarray*}
E_{o,o} [f(\theta,l_1,X,\tilde{X}) ]&\leq&\sum_{k=0}^\infty
K_1^{1/2}\bigl(1+k\ee^{a_1|\theta|k}\bigr)^{1/2}K_2^{1/2}\ee^{-a_3k/2}\\
&\leq&2(K_1K_2)^{1/2}\sum_{k=0}^\infty k^{1/2}\ee^{(a_1|\theta
|-a_3)k/2}\\
&<&\infty,
\end{eqnarray*}
provided that $|\theta|<a_3/a_1$.

The constant $c_2$ is chosen in \cite{YilmazQequalsA} such that it
satisfies $c_2<\min(c_1/4, a_3/a_1)$, along with a few other
conditions. Thus, (\ref{lemsamikibin}) holds for every $\theta\in
\mathcal{C}_{\eq}=\mathcal{C}(c_2)$.
\end{pf}
%

\subsection{\texorpdfstring{Proof of Theorem \protect\ref{existence}}{Proof of Theorem 3.1}}

For every $n\geq2$, $\theta\in\mathcal{C}_{\eq}$ and $\omega\in\Omega$,
%
%
\begin{eqnarray}\label{basel}\qquad
h_n(\theta,\omega)&=&E_o^\omega[\exp\{\langle\theta
,X_{H_n}\rangle-\Lambda_a(\theta)H_n \},H_n=\tau_k \mbox{ for some
}k\geq1 ]\nonumber\\
&=&\sum_{z\in U}E_o^\omega[\exp\{\langle\theta,X_{H_n}\rangle
-\Lambda_a(\theta)H_n \},\nonumber\\
&&\hspace*{36pt}X_1=z,H_n=\tau_k \mbox{ for some }k\geq
1 ]\nonumber\\
&=&\sum_{z\in U} \pi(0,z) \exp\{-\Lambda_a(\theta)\}\nonumber\\
&&\hspace*{14pt}{}\times E_z^\omega[\exp
\{\langle\theta,X_{H_n}\rangle-\Lambda_a(\theta)H_n \}
,H_n=\tau_k \mbox{ for some }k\geq1 ]\nonumber\\
&=&\sum_{z\in U}\pi(0,z)\exp\{\langle\theta,z\rangle-\Lambda_a(\theta)\}
h_{n-\langle z,e_1\rangle}(\theta,T_z\omega).
\end{eqnarray}
Here, (\ref{basel}) is obtained by shifting the environment by $z$.

Define a new function $\bar{h}_n(\theta,\cdot)\dvtx\Omega\to\mathbb{R}$ by
%
%
\begin{equation}\label{analojiyap}
\bar{h}_n(\theta,\omega):=\frac{1}{n-1}\sum_{i=2}^nh_i(\theta,\omega).
\end{equation}
Since $(\bar{h}_n(\theta,\cdot))_{n\geq1}$ is bounded in $L^2(\mathbb
{P})$ by (\ref{lemsamikibin}), it has a subsequence $(\bar
{h}_{n_k}(\theta,\cdot))_{k\geq1}$ that converges weakly to some
$h(\theta,\cdot)\in L^2(\mathbb{P})$.

It follows immediately from (\ref{basel}) that
%
%
\begin{eqnarray}\label{acesko}
\bar{h}_n(\theta,\omega)&=&\sum_{z\in U}\pi(0,z)\exp\{\langle\theta
,z\rangle-\Lambda_a(\theta)\}\bar{h}_n(\theta,T_z\omega)\nonumber\\
&&{} +\frac{1}{n-1} \bigl(h_1(\theta,T_{e_1}\omega)-h_2(\theta
,T_{-e_1}\omega)\\
&&\hspace*{41pt}{} -h_n(\theta,T_{e_1}\omega)+h_{n+1}(\theta,T_{-e_1}\omega
) \bigr).\nonumber
\end{eqnarray}
Set $n=n_k$ and take the weak limit of both sides of (\ref{acesko}) as
$k\to\infty$. Since the term on the second line converges (strongly
and, hence, weakly) to zero in $L^2(\mathbb{P})$, we conclude that
%
%
\begin{equation}\label{isitisit}\qquad
h(\theta,\omega)=\sum_{z\in U}\pi(0,z)\exp\{\langle\theta,z\rangle
-\Lambda_a(\theta)\}h(\theta,T_z\omega) \qquad\mbox{for $\mathbb
{P}$-a.e. $\omega$.}
\end{equation}

Note that $h(\theta,\omega)\geq0$ for $\mathbb{P}$-a.e. $\omega$.
Equation (\ref{isitisit}) [in combination with (\ref{ellipticity})]
implies that the set $\{\omega\in\Omega\dvtx h(\theta,\omega)=0\}$ is
invariant under $(T_z)_{z\in U}$. Since (\ref{i.i.d.}) ensures that the
environment is ergodic under these shifts, $\mathbb{P}(h(\theta,\cdot
)=0)\in\{0,1\}$. However, $\mathbb{E}\{h(\theta,\cdot)\}>0$ by (\ref
{lemsambin}). Therefore, we conclude that \mbox{$\mathbb{P}(h(\theta,\cdot
)>0)=1$}. We have thus proven Theorem \ref{existence}.
%

\subsection{A useful representation}

Define a function $\varphi\dvtx\mathcal{C}_{\eq}\times\Omega\times
\mathbb
{Z}^d\to\mathbb{R}^+$ by setting
%
%
\begin{equation}\label{fifi}
\varphi(\theta,\omega,x):=\frac{E_o^\omega[\exp\{\langle\theta
,X_{\tau_1}\rangle-\Lambda_a(\theta)\tau_1 \},X_{\tau_1}=x
]}{P_o^{T_x\omega}(\beta=\infty)}
\end{equation}
for every $\theta\in\mathcal{C}_{\eq}$, $\omega\in\Omega$ and $x\in
\mathbb{Z}^d$. Note that $\varphi(\theta,\omega,x)=0$ unless $\langle
x,e_1\rangle\geq1$.

The following lemma will be useful in the next section.
\begin{lemma}\label{kursunlem}
For every $\theta\in\mathcal{C}_{\eq}$, there exists a $\mathcal
{B}_o^+(e_1)$-measurable $g(\theta,\cdot)\in L^2(\mathbb{P})$ such that
%
%
\begin{equation}\label{kursunlar}
h(\theta,\omega)=\sum_{x\in\mathbb{Z}^d}\varphi(\theta,\omega,x)g(\theta
,T_x\omega) \qquad\mbox{for $\mathbb{P}$-a.e. $\omega$.}
\end{equation}
\end{lemma}
\begin{pf}
Recall (\ref{kafes}) and (\ref{kameriye}). For every $n\geq2$ and
$\theta\in\mathcal{C}_{\eq}$, define $\bar{g}_n(\theta,\cdot)\in
L^2(\mathbb{P})$ analogously to (\ref{analojiyap}). Since $(\bar
{g}_n(\theta,\cdot))_{n\geq1}$ is bounded in $L^2(\mathbb{P})$ by (\ref
{lemsamiki}), it has a subsequence $(\bar{g}_{n_k}(\theta,\cdot))_{k\geq
1}$ that converges weakly to some $g(\theta,\cdot)\in L^2(\mathbb{P})$.
[Choose $(n_k)_{k\geq1}$ to be a further subsequence of the subsequence
in the proof of Theorem \ref{existence} so that $(\bar{h}_{n_k}(\theta
,\cdot))_{k\geq1}$ converges weakly to $h(\theta,\cdot)\in L^2(\mathbb
{P})$.] Note that $g_n(\theta,\cdot)$ is $\mathcal
{B}_o^+(e_1)$-measurable for every $n\geq1$ since the event $\{\beta
=\infty\}$ is part of the definition of $g_n(\theta,\cdot)$. Hence,
$g(\theta,\cdot)$ is $\mathcal{B}_o^+(e_1)$-measurable.

For every $N\geq1$, $n\geq N$, $\theta\in\mathcal{C}_{\eq}$ and $\omega
\in\Omega$,
%
%
\begin{eqnarray}\label{dayangenc}\quad
h_n(\theta,\omega)&=&E_o^\omega[\exp\{\langle\theta
,X_{H_n}\rangle-\Lambda_a(\theta)H_n \},H_n=\tau_k \mbox{ for some
}k\geq1 ]\nonumber\\
&=&E_o^\omega[\exp\{\langle\theta,X_{H_n}\rangle-\Lambda
_a(\theta)H_n \},\nonumber\\
&&\hspace*{17.8pt}|X_{\tau_1}|\geq N,H_n=\tau_k \mbox{ for some
}k\geq1 ]\\
&&{} +\sum_{|x|<N}E_o^\omega[\exp\{\langle\theta
,X_{H_n}\rangle-\Lambda_a(\theta)H_n \},\nonumber\\
&&\hspace*{56.2pt}X_{\tau_1}=x,H_n=\tau_k
\mbox{ for some }k\geq1 ]\nonumber.
\end{eqnarray}
Denote the first term in (\ref{dayangenc}) by $R_{N,n}(\theta,\omega)$.
It follows immediately from (\ref{yalebbim}), the renewal structure and
the monotone convergence theorem that
%
%
\begin{eqnarray}\label{busizi}
&&
\lim_{N\to\infty}\sup_{n\geq N}\mathbb{E}\{R_{N,n}(\theta,\cdot)\}\nonumber\\[-8pt]\\[-8pt]
&&\qquad\leq
\lim_{N\to\infty}E_o [\exp\{\langle\theta,X_{\tau_1}\rangle
-\Lambda_a(\theta)\tau_1 \},|X_{\tau_1}|\geq N ]=0.\nonumber
\end{eqnarray}
Recall (\ref{fifi}) and observe that, for every $n\geq N$,
\begin{eqnarray*}
h_n(\theta,\omega)&=&R_{N,n}(\theta,\omega)\\
&&{}+\sum_{|x|<N}E_o^\omega
[\exp\{\langle\theta,X_{H_n}\rangle-\Lambda_a(\theta)H_n \},\\
&&\hspace*{56.2pt}X_{\tau_1}=x,H_n=\tau_k \mbox{ for some }k\geq1 ]\\
&=&R_{N,n}(\theta,\omega)\\
&&{}+\sum_{|x|<N}E_o^\omega[\exp\{\langle
\theta,X_{\tau_1}\rangle-\Lambda_a(\theta)\tau_1 \},X_{\tau
_1}=x ]\\
&&\hspace*{35.3pt}{}\times\ee^{-\langle\theta,x\rangle}E_x^\omega
[\exp\{\langle\theta,X_{H_n}\rangle-\Lambda_a(\theta
)H_n \},\\
&&\hspace*{97.3pt}H_n=\tau_k \mbox{ for some }k\geq1 |\beta=\infty
]\\
&=&R_{N,n}(\theta,\omega)+\sum_{|x|<N}\varphi(\theta,\omega
,x)g_{n-\langle x,e_1\rangle}(\theta,T_x\omega).
\end{eqnarray*}
Therefore, whenever $n_k\geq N$,
%
%
\begin{eqnarray}\label{sorin}
\frac{1}{n_k}\sum_{i=N}^{n_k}h_i(\theta,\omega)&=&\frac{1}{n_k}\sum
_{i=N}^{n_k}R_{N,i}(\theta,\omega)\nonumber\\[-8pt]\\[-8pt]
&&{}+\sum_{|x|<N}\varphi(\theta,\omega
,x)\frac{1}{n_k}\sum_{i=N}^{n_k}g_{i-\langle x,e_1\rangle}(\theta
,T_x\omega).\nonumber
\end{eqnarray}
Multiplying both sides of (\ref{sorin}) by any indicator function $\chi
\in L^\infty(\mathbb{P})$, integrating against $\mathbb{P}$ and letting
$k$ tend to infinity, we arrive at the following inequality:
\[
\biggl|\int h(\theta,\omega)\chi(\omega)\,\dd\mathbb{P} - \int\sum
_{|x|<N}\varphi(\theta,\omega,x)g(\theta,T_x\omega)\chi(\omega)\,\dd
\mathbb{P} \biggr|\leq\sup_{n\geq N}\mathbb{E}\{R_{N,n}(\theta,\cdot
)\}.
\]
Finally, let $N$ tend to infinity. The monotone convergence theorem and
(\ref{busizi}) imply the desired result.
\end{pf}
%

\section{Proof of our results on averaged large deviations}\label{kangiiki}

We will start this section by stating two results concerning the unique
minimizer of Varadhan's variational formula (\ref{divaneasik}). We will
then give a series of lemmas. Finally, we will combine everything and
prove Theorem \ref{varminvar}.

%

\subsection{The unique minimizer of Varadhan's variational
formula}\label{evinara}

Assume (\ref{kondisin}). Take any $\xi\in\mathcal{A}_a$. Since the
Hessian $\mathcal{H}_a$ of $\Lambda_a$ is positive definite on $\mathcal
{C}_a$ by Theorem \ref{Caveraged}, there exists a unique $\theta\in
\mathcal{C}_a$ satisfying $\xi=\nabla\Lambda_a(\theta)$. In the next
paragraph, we define a probability measure $\hat{\mu}_\xi^\infty\in
M_1(\Omega\times U^{\mathbb{N}})$ by specifying the integrals of
certain test functions against this measure.

For every $N,M,K\geq0$, take any bounded function $f\dvtx\Omega\times
U^{\mathbb{N}}\rightarrow\mathbb{R}$ such that $f(\cdot,(z_i)_{i\geq
1})$ is independent of $(z_i)_{i>K}$ and is measurable with respect to
%
%
\begin{equation}\label{kuruyorum}
\mathcal{B}_N^M(e_1)=\mathcal{B}_N^M:=\sigma(\omega_x\dvtx-N\leq\langle
x,e_1\rangle\leq M).
\end{equation}
Define $\hat{\mu}_\xi^\infty\in M_1(\Omega\times U^{\mathbb{N}})$ by setting
%
%
\begin{eqnarray}\label{muyucananan}
\int f\,\dd\hat{\mu}_\xi^\infty&:=&\sum_{j=0}^\infty
E_o [\tau_N\leq j<\tau_{N+1}, \nonumber\\[-4pt]
&&\hspace*{33.4pt}f(T_{X_j}\omega,Z_{j+1}^\infty
)\exp\{\langle\theta,X_{\tau_J}\rangle- \Lambda_a(\theta)\tau_J\}
|\beta=\infty]\hspace*{-22pt}\\
&&{}\times \bigl(E_o [\tau_1\exp\{\langle\theta
,X_{\tau_1}\rangle- \Lambda_a(\theta)\tau_1\} |\beta=\infty]\bigr)^{-1},\nonumber
\end{eqnarray}
where $J:=N+M+K+1$ and $Z_{j+1}^\infty=(Z_{j+i})_{i\geq
1}:=(X_i-X_{i-1})_{i\geq1}$. (The measure $\hat{\mu}_\xi^\infty$ is
well defined. See the proof of Theorem \ref{averagedconditioningsh}.)

The following theorem states that the empirical process
\[
\nu_{n,X}^\infty:= \frac{1}{n}\sum_{k=0}^{n-1}\one_{T_{X_k}\omega
,Z_{k+1}^\infty}
\]
of the walk under $P_o$ converges to $\hat{\mu}_\xi^\infty$ when the
particle is conditioned to have mean velocity $\xi$.
\begin{theorem}\label{averagedconditioningsh}
Assume (\ref{kondisin}). For every $\xi\in\mathcal{A}_a$, $\varepsilon>0$,
$N,M,K\geq0$ and $f\dvtx\Omega\times U^{\mathbb{N}}\rightarrow\mathbb{R}$
bounded such that $f(\cdot,(z_i)_{i\geq1})$ is independent of
$(z_i)_{i>K}$ and is $\mathcal{B}_N^M$-measurable, the following holds:
%
%
\begin{eqnarray}\label{saglik}
&&\limsup_{\delta'\to0}\limsup_{n\rightarrow\infty}\frac{1}{n}\log
P_o \biggl(\biggl|\int f\,\dd\nu_{n,X}^\infty-\int f\,\dd\hat{\mu}_\xi^\infty
\biggr|>\varepsilon \Big|\nonumber\\[-8pt]\\[-8pt]
&&\hspace*{138.4pt}\hspace*{40.6pt} \biggl|\frac{X_n}{n}-\xi\biggr|\leq
\delta' \biggr)<0.\nonumber
\end{eqnarray}
\end{theorem}
\begin{pf}
Definition 9 of \cite{YilmazAveraged} introduces a probability measure
$\bar{\mu}_\xi^\infty\in M_1(U^{\mathbb{N}})$ by the formula in (\ref
{muyucananan}), except that the test functions do not depend on $\omega
$ in that case. Proposition 16 of \cite{YilmazAveraged} shows that $\bar
{\mu}_\xi^\infty$ is well defined, and Theorem 17 of \cite
{YilmazAveraged} establishes the analog of (\ref{saglik}) for $\bar{\mu
}_\xi^\infty$. The proofs of these results generalize to our setting
without any nontrivial change. Therefore, we omit the proof of
Theorem~\ref{averagedconditioningsh}. Also, note that Theorem \ref
{averagedconditioningsh} is proved in \cite{YilmazSpaceTime} for the
related model of space--time RWRE.
\end{pf}

In the first paragraph of Section \ref{unminvar}, we mentioned an
existence and uniqueness result for the minimizer of Varadhan's
variational formula (\ref{divaneasik}) for the averaged rate function
$I_a$. The following is the precise statement.
\begin{theorem}\label{cizme}
Assume (\ref{kondisin}). For every $\xi\in\mathcal{A}_a$, $\hat{\mu}_\xi
^\infty$ induces a $\mathbb{Z}^d$-valued transient process with
stationary and ergodic increments in $U$ via the map
\[
(\omega,z_1,z_2,z_3,\ldots)\mapsto(z_1,z_1+z_2,z_1+z_2+z_3,\ldots).
\]
Extend this process to a probability measure on doubly infinite paths
$(x_i)_{i\in\mathbb{Z}}$ and refer to its restriction to $W_{\infty
}^{\mathrm{tr}}$ as $\mu_\xi^\infty$. With this notation, $\mu_\xi^\infty$ is
the unique minimizer of Varadhan's variational formula (\ref{divaneasik}).
\end{theorem}
\begin{pf}
This is Theorem 10 of \cite{YilmazAveraged}, with the following
difference: that result is concerned with $\bar{\mu}_\xi^\infty$ (which
was mentioned in the proof of Theorem \ref{averagedconditioningsh}) and
it uses the map
\[
(z_1,z_2,z_3,\ldots)\mapsto(z_1,z_1+z_2,z_1+z_2+z_3,\ldots)
\]
to induce a $\mathbb{Z}^d$-valued transient process with stationary and
ergodic increments in $U$. However, since $\bar{\mu}_\xi^\infty$ is the
marginal of $\hat{\mu}_\xi^\infty$ on $U^{\mathbb{N}}$, the $\mathbb
{Z}^d$-valued process induced by $\bar{\mu}_\xi^\infty$ is nothing but
$\mu_\xi^\infty$.
\end{pf}
%

\subsection{The Markovian structure of the minimizer}\label{evinye}

Assume (\ref{kondisin}). Take any $\xi\in\mathcal{A}_{\eq}$. Let $\theta
\in\mathcal{C}_{\eq}$ denote the unique solution of $\xi=\nabla\Lambda
_a(\theta)$. Recall the environment kernel $\hat{\pi}^\theta$ defined
in (\ref{azdinyn}) via h-transform. For any $x\in\mathbb{Z}^d$ and
$\omega\in\Omega$, abbreviate the notation introduced in Definition \ref
{ortamkeli} by writing $P_x^{\theta,\omega}$ and $E_x^{\theta,\omega}$
instead of $P_x^{\hat{\pi}^\theta,\omega}$ and $E_x^{\hat{\pi}^\theta
,\omega}$, respectively.

For every $n\geq1$, define $\hat{\mu}_{n,\xi}^\infty\in M_1(\Omega
\times U^{\mathbb{N}})$ as follows:
%
%
\begin{equation}\label{gurbetlik}
\hat{\mu}_{n,\xi}^\infty(\cdot):=\frac{\mathbb{E}\{h(\theta,\omega
)P_o^{\theta,\omega}((T_{X_n}\omega,Z_{n+1}^\infty)\in\cdot)\}}{\mathbb
{E}\{h(\theta,\omega)\}}.
\end{equation}
\begin{lemma}
For every $N,M,K\geq0$ and $f\dvtx\Omega\times U^{\mathbb
{N}}\rightarrow
\mathbb{R}$ bounded such that $f(\cdot,(z_i)_{i\geq1})$ is independent
of $(z_i)_{i>K}$ and is $\mathcal{B}_N^M$-measurable,
%
%
\begin{equation}\label{komurhan}
\int f\,\dd\hat{\mu}_{n,\xi}^\infty=\frac{E_o[f(T_{X_n}\omega
,Z_{n+1}^\infty)\exp\{\langle\theta,X_{\tau_L}\rangle-\Lambda_a(\theta
)\tau_L\}]}{E_o[\exp\{\langle\theta,X_{\tau_1}\rangle-\Lambda_a(\theta
)\tau_1\}]}
\end{equation}
for every $L\geq n+M+K+1$.
\end{lemma}
\begin{pf}
For every $N,M,K\geq0$, take a bounded function $f\dvtx\Omega\times
U^{\mathbb{N}}\rightarrow\mathbb{R}$ such that $f(\cdot,(z_i)_{i\geq
1})$ is independent of $(z_i)_{i>K}$ and is $\mathcal
{B}_N^M$-measurable; see (\ref{kuruyorum}). For every $L\geq n+M+K+1$,
%
%
\begin{eqnarray}\qquad
&&\mathbb{E}\{h(\theta,\omega)\}\int f\,\dd\hat{\mu}_{n,\xi
}^\infty\nonumber\\
&&\qquad=\mathbb{E}\{h(\theta,\omega)E_o^{\theta,\omega}[f(T_{X_n}\omega
,Z_{n+1}^\infty)]\}\nonumber\\
\label{garbir}
&&\qquad =\mathbb{E} \biggl\{h(\theta,\omega)E_o^\omega\biggl[f(T_{X_n}\omega
,Z_{n+1}^\infty)\nonumber\\[-8pt]\\[-8pt]
&&\qquad\quad\hspace*{63pt}{}\times\exp\{\langle\theta,X_{n+K}\rangle-\Lambda_a(\theta
)(n+K)\}\frac{h(\theta,T_{X_{n+K}}\omega)}{h(\theta,\omega)}
\biggr] \biggr\}\nonumber\\
&&\qquad =E_o[f(T_{X_n}\omega,Z_{n+1}^\infty)\exp\{\langle\theta
,X_{n+K}\rangle-\Lambda_a(\theta)(n+K)\}h(\theta,T_{X_{n+K}}\omega
)]\nonumber\\
\label{garuc}
&&\qquad =E_o[f(T_{X_n}\omega,Z_{n+1}^\infty)\exp\{\langle\theta
,X_{H_L}\rangle-\Lambda_a(\theta)H_L\}h(\theta,T_{X_{H_L}}\omega)]\\
\label{gardort}
&&\qquad =\sum_{\langle x,e_1\rangle\geq1}E_o[f(T_{X_n}\omega
,Z_{n+1}^\infty)\exp\{\langle\theta,X_{H_L}\rangle-\Lambda_a(\theta
)H_L\}\nonumber\\[-8pt]\\[-8pt]
&&\qquad\quad\hspace*{90.4pt}{}\times\varphi(\theta,T_{X_{H_L}}\omega,x)g(\theta,T_{X_{H_L}+x}\omega)].\nonumber
\end{eqnarray}
\textit{Explanation}: (\ref{garbir}) follows from the definition of $\hat{\pi
}^\theta$ by noting that $f(T_{X_n}\omega,Z_{n+1}^\infty)$ depends only
on the first $n+K$ steps of the walk. (\ref{garuc}) holds because $H_L$
is a stopping time and $H_L\geq n+K$. The representation of $h(\theta
,\cdot)$ in (\ref{kursunlar}) gives (\ref{gardort}).

For any $x\in\mathbb{Z}^d$ such that $\langle x,e_1\rangle\geq1$,
%
%
\begin{eqnarray}
&&
E_o[f(T_{X_n}\omega,Z_{n+1}^\infty)\nonumber\\
&&\hspace*{2.7pt}\quad{}\times\exp\{\langle\theta,X_{H_L}\rangle
-\Lambda_a(\theta)H_L\}\varphi(\theta,T_{X_{H_L}}\omega,x)g(\theta
,T_{X_{H_L}+x}\omega)]\nonumber\\
\label{sonbir}
&&\qquad =\sum_{\langle y,e_1\rangle=L} \mathbb{E}\bigl\{E_o^\omega
[f(T_{X_n}\omega,Z_{n+1}^\infty)\exp\{\langle\theta,y\rangle-\Lambda
_a(\theta)H_L\}\nonumber\\
&&\qquad\quad\hspace*{126.4pt}{}\times\varphi(\theta,T_y\omega,x),X_{H_L}=y]\\
&&\qquad\quad\hspace*{126.4pt}\hspace*{49.3pt}{}\times g(\theta
,T_{y+x}\omega)\bigr\}\nonumber\\
\label{soniki}
&&\qquad =\sum_{\langle y,e_1\rangle=L} E_o[f(T_{X_n}\omega
,Z_{n+1}^\infty)\exp\{\langle\theta,y\rangle-\Lambda_a(\theta)H_L\}\nonumber\\[-8pt]\\[-8pt]
&&\qquad\quad\hspace*{112.4pt}{}\times
\varphi(\theta,T_y\omega,x),X_{H_L}=y]\mathbb{E}\{g(\theta,\omega)\}
\nonumber\\
\label{sonuc}
&&\qquad =E_o[f(T_{X_n}\omega,Z_{n+1}^\infty)\exp\{\langle\theta
,X_{H_L}\rangle-\Lambda_a(\theta)H_L\}\varphi(\theta,T_{X_{H_L}}\omega
,x)]\nonumber\\[-8pt]\\[-8pt]
&&\qquad\quad\hspace*{0pt}{}\times\mathbb{E}\{g(\theta,\omega)\}.\nonumber
\end{eqnarray}
\textit{Explanation}: for any $y\in\mathbb{Z}^d$ such that $\langle y,e_1\rangle
=L$, the random quantities $E_o^\omega[ \ldots,X_{H_L}=y]$ and
$g(\theta,T_{y+x}\omega)$ appearing in (\ref{sonbir}) are independent
because the former is\vadjust{\goodbreak}
measurable with respect to $\sigma(\omega_{x'}\dvtx
\langle x'-x,e_1\rangle<L)$, whereas the latter is $\mathcal
{B}_{L+\langle x,e_1\rangle}^+(e_1)$-measurable; see Lemma \ref
{kursunlem}. This independence (in combination with the stationarity of
$\mathbb{P}$) gives (\ref{soniki}).

By plugging (\ref{sonuc}) into (\ref{gardort}), we see that
%
%
\begin{eqnarray}\qquad
&&\frac{\mathbb{E}\{h(\theta,\cdot)\}}{\mathbb{E}\{g(\theta,\cdot)\}}\int
f\,\dd\hat{\mu}_{n,\xi}^\infty\nonumber\\
\label{leybilbir}
&&\qquad =\sum_{\langle x,e_1\rangle\geq1}\sum_{k=1}^{L+1}\mathbb{E}
\bigl\{E_o^\omega[f(T_{X_n}\omega,Z_{n+1}^\infty)\exp\{\langle\theta
,X_{\tau_k}\rangle-\Lambda_a(\theta)\tau_k\},\nonumber\\
&&\qquad\quad\hspace*{157pt}\hspace*{-3.8pt}\tau_{k-1}\leq H_L<\tau
_k,X_{\tau_k}=x]\\
&&\qquad\quad\hspace*{21.6pt}\hspace*{151pt}{}\times\bigl({P_o^{T_x\omega}(\beta=\infty)} \bigr)^{-1}\bigr\}\nonumber\\
\label{leybiliki}
&&\qquad =\frac{1}{P_o(\beta=\infty)}\nonumber\\
&&\qquad\quad{}\times\sum_{\langle x,e_1\rangle\geq1}\sum
_{k=1}^{L+1}E_o[f(T_{X_n}\omega,Z_{n+1}^\infty)\exp\{\langle\theta
,X_{\tau_k}\rangle-\Lambda_a(\theta)\tau_k\},\\
&&\qquad\quad\hspace*{153pt}\tau_{k-1}\leq H_L<\tau
_k,X_{\tau_k}=x]\nonumber\\
\label{leybiluc}
&&\qquad =\frac{1}{P_o(\beta=\infty)}E_o[f(T_{X_n}\omega,Z_{n+1}^\infty
)\exp\{\langle\theta,X_{\tau_L}\rangle-\Lambda_a(\theta)\tau_L\}].
\end{eqnarray}
\textit{Explanation}: (\ref{leybilbir}) follows from the definition of $\varphi
(\theta,\cdot,\cdot)$ given in (\ref{fifi}). For every $x\in\mathbb
{Z}^d$ such that $\langle x,e_1\rangle\geq1$, the random quantity
$P_o^{T_x\omega}(\beta=\infty)$ is independent of the ratio $\frac
{E_o^\omega[ \ldots,X_{\tau_k}=x]}{P_o^{T_x\omega}(\beta=\infty)}$
appearing in (\ref{leybilbir}) since the latter is easily seen to be
equal to an expectation involving the stopping time $H_{\langle
x,e_1\rangle}$ (and nothing beyond that). This independence implies
(\ref{leybiliki}). Using (\ref{yalebbim}), the $\tau_k$ in the
exponential can be replaced first by $\tau_{L+1}$ and then by $\tau
_{L}$. This gives (\ref{leybiluc}).

Finally, observe that (\ref{leybiluc}) agrees with (\ref{komurhan}),
except that the normalization constant has to be simplified. However,
it is clear that the constant in (\ref{komurhan}) is correct [take
$f\equiv1$ and apply (\ref{yalebbim})].
\end{pf}
\begin{lemma}\label{zafiyet}
For every $f\dvtx\Omega\times U^{\mathbb{N}}\rightarrow\mathbb{R}$ bounded
such that $f(\cdot,(z_i)_{i\geq1})$ is $\mathcal{B}$-measurable, the
following convergence takes place:
\[
\lim_{n\to\infty}\int f\,\dd\hat{\mu}_{n,\xi}^\infty=\int
f\,\dd\hat{\mu}_\xi^\infty.
\]
In particular, $(\hat{\mu}_{n,\xi}^\infty)_{n\geq1}$ converges weakly
to $\hat{\mu}_\xi^\infty$.
\end{lemma}
\begin{pf}
For any $N,M,K\geq0$, take a bounded function $f\dvtx\Omega\times
U^{\mathbb
{N}}\rightarrow\mathbb{R}$ such that $f(\cdot,(z_i)_{i\geq1})$ is
independent of $(z_i)_{i>K}$ and $\mathcal{B}_N^M$-measurable. Let $J:=N+M+K+1$.
%
%
\begin{eqnarray}\label{atasgg}\hspace*{33pt}
&&E_o[\exp\{\langle\theta,X_{\tau_1}\rangle-\Lambda_a(\theta)\tau_1\}
]\lim_{n\to\infty}\int f\,\dd\hat{\mu}_{n,\xi}^\infty\\
\label{ozgurolcam}
&&\qquad =\lim_{n\to\infty}E_o[n<\tau_N, f(T_{X_n}\omega,Z_{n+1}^\infty
)\exp\{\langle\theta,X_{\tau_J}\rangle-\Lambda_a(\theta)\tau_J\}]\\
&&\qquad\quad{} +\lim_{n\to\infty}\sum_{i=0}^\infty E_o[\tau_{N+i}\leq
n<\tau_{N+i+1},\nonumber\\
&&\qquad\hspace*{80.2pt}f(T_{X_n}\omega,Z_{n+1}^\infty)\exp\{\langle\theta
,X_{\tau_{J+i}}\rangle-\Lambda_a(\theta)\tau_{J+i}\}]\nonumber\\
\label{danso}
&&\qquad =\lim_{n\to\infty}\sum_{j=0}^n \Biggl(\sum_{i=0}^\infty E_o[\tau
_i=n-j,\exp\{\langle\theta,X_{\tau_i}\rangle-\Lambda_a(\theta)\tau_i\}
] \Biggr)\\
&&\qquad\quad\hspace*{37.3pt}{}\times E_o[\tau_N\leq j<\tau_{N+1},\nonumber\\
&&\qquad\hspace*{78.1pt}f(T_{X_j}\omega
,Z_{j+1}^\infty)\exp\{\langle\theta,X_{\tau_J}\rangle-\Lambda_a(\theta
)\tau_J\} |\beta=\infty]\nonumber\\
&&\qquad =S(\theta)\sum_{j=0}^\infty E_o[\tau_N\leq j<\tau
_{N+1},\nonumber\\
&&\qquad\hspace*{65.1pt}f(T_{X_j}\omega,Z_{j+1}^\infty)\exp\{\langle\theta,X_{\tau
_J}\rangle-\Lambda_a(\theta)\tau_J\} |\beta=\infty]\nonumber\\
\label{atesgg}
&&\qquad =S(\theta)E_o [\tau_1\exp\{\langle\theta,X_{\tau
_1}\rangle- \Lambda_a(\theta)\tau_1\} |\beta=\infty]\int
f\,\dd\hat{\mu}_\xi^\infty.
\end{eqnarray}
\textit{Explanation}: (\ref{ozgurolcam}) follows from (\ref{komurhan}). The
first term in (\ref{ozgurolcam}) goes to zero as $n\to\infty$ by the
dominated convergence theorem. The renewal theorem for aperiodic
sequences (see \cite{BreimanBook}, Theorem 10.8) implies that the sum
$(\sum_{i=0}^\infty E_o[\cdots])$ in (\ref{danso}) converges to some
constant $S(\theta)$ as $n\to\infty$. Observe that the constants in
(\ref{atasgg}) and (\ref{atesgg}) have to agree because $\hat{\mu
}_{n,\xi}^\infty$ and $\hat{\mu}_\xi^\infty$ are known to be
probability measures.

Thus far, we have shown that $\lim_{n\to\infty}\int f\,\dd\hat
{\mu}_{n,\xi}^\infty=\int f\,\dd\hat{\mu}_\xi^\infty$ for a
separating class of test functions. However, this is sufficient to
conclude that $(\hat{\mu}_{n,\xi}^\infty)_{n\geq1}$ converges weakly to
$\hat{\mu}_\xi^\infty$ since $M_1(\Omega\times U^\mathbb{N})$ is compact.
\end{pf}

Let $\mathbb{Q}_\xi\in M_1(\Omega)$ be the marginal of $\hat{\mu}_\xi
^\infty\in M_1(\Omega\times U^\mathbb{N})$.
\begin{lemma}\label{kiliyir}
$\mathbb{Q}_\xi$ is $\hat{\pi}^\theta$-invariant, that is,
\[
\sum_{z\in U}\dd\mathbb{Q}_\xi(T_{-z}\omega)\hat{\pi}^\theta
(T_{-z}\omega,z)=\dd\mathbb{Q}_\xi(\omega).
\]
\end{lemma}
\begin{pf}
For any $f\in L^\infty(\mathbb{P})$, define $\hat{\pi}^\theta f\dvtx
\Omega
\to\mathbb{R}$ in the usual way:
\[
(\hat{\pi}^\theta f)(\omega):=\sum_{z\in U}\hat{\pi}^\theta(\omega
,z)f(T_z\omega).
\]
Recall (\ref{gurbetlik}). For every $n\geq1$,
\begin{eqnarray*}
\int(\hat{\pi}^\theta f)\,\dd\hat{\mu}_{n,\xi}^\infty&=&\frac{\mathbb
{E}\{h(\theta,\omega)E_o^{\theta,\omega}[(\hat{\pi}^\theta
f)(T_{X_n}\omega)]\}}{\mathbb{E}\{h(\theta,\omega)\}}\\
&=&\frac{\mathbb{E}\{
h(\theta,\omega)E_o^{\theta,\omega}[f(T_{X_{n+1}}\omega)]\}}{\mathbb
{E}\{h(\theta,\omega)\}}\\
&=&\int f\,\dd\hat{\mu}_{n+1,\xi}^\infty,
\end{eqnarray*}
by the Markov property. Let $n$ tend to infinity, use Lemma \ref
{zafiyet} and conclude that
\[
\int(\hat{\pi}^\theta f)\,\dd\mathbb{Q}_\xi=\int(\hat{\pi}^\theta
f)\,\dd\hat{\mu}_\xi^\infty=\int f\,\dd\hat{\mu}_\xi^\infty
=\int f\,\dd\mathbb{Q}_\xi.
\]
This is equivalent to the desired result.
\end{pf}
\begin{lemma}\label{markovkopek}
$\hat{\mu}_\xi^\infty$ induces, via the map
\[
(\omega, z_1,z_2,z_3,\ldots)\mapsto(\omega,T_{z_1}\omega
,T_{z_1+z_2}\omega,T_{z_1+z_2+z_3}\omega,\ldots),
\]
an $\Omega$-valued stationary Markov process with marginal $\mathbb
{Q}_\xi$ and transition kernel
\[
\overline{\pi}^\theta(\omega,\omega'):=\sum_{z\dvtx T_z\omega=\omega
'}\hat
{\pi}^\theta(\omega,z).
\]
\end{lemma}
\begin{pf}
For any $n\geq1$, $K\geq0$ and any two bounded measurable functions
$f\dvtx\Omega^{K+1}\rightarrow\mathbb{R}$ and $g\dvtx\Omega^{\mathbb
{N}}\rightarrow\mathbb{R}$, it follows from (\ref{gurbetlik}) and the
Markov property that
\begin{eqnarray*}
&&\mathbb{E}\{h(\theta,\omega)\}\int f(\omega,T_{z_1}\omega,\ldots
,T_{z_1+\cdots+z_K}\omega)\\
&&\qquad\quad\hspace*{25.7pt}{}\times g(T_{z_1+\cdots+z_K}\omega,T_{z_1+\cdots
+z_{K+1}}\omega,\ldots)\,\dd\hat{\mu}_{n,\xi}^\infty\\
&&\qquad =\mathbb{E}\{h(\theta,\omega)E_o^{\theta,\omega}[f(T_{X_n}\omega
,\ldots,T_{X_{n+K}}\omega)g(T_{X_{n+K}}\omega,T_{X_{n+K+1}}\omega,\ldots
)]\}\\
&&\qquad =\mathbb{E} \bigl\{h(\theta,\omega)E_o^{\theta,\omega}
\bigl[f(T_{X_n}\omega,\ldots,T_{X_{n+K}}\omega)\\
&&\qquad\quad\hspace*{67.5pt}{}\times E_o^{\theta,\omega
}[g(T_{X_{n+K}}\omega,T_{X_{n+K+1}}\omega,\ldots)|X_{n+K}] \bigr]
\bigr\}\\
&&\qquad =\mathbb{E}\{h(\theta,\omega)\}\\
&&\qquad\quad{}\times\int f(\omega,T_{z_1}\omega,\ldots
,T_{z_1+\cdots+z_K}\omega)\\
&&\qquad\quad\hspace*{20.9pt}{}\times E_o^{\theta,T_{z_1+\cdots+z_K}\omega
}[g(T_{z_1+\cdots+z_K}\omega,T_{z_1+\cdots+z_K+X_1}\omega,\ldots
)]\,\dd\hat{\mu}_{n,\xi}^\infty.
\end{eqnarray*}
Let $n$ tend to infinity, use Lemma \ref{zafiyet} and conclude that
$\hat{\mu}_\xi^\infty$ indeed induces an $\Omega$-valued Markov process
with marginal $\mathbb{Q}_\xi$ and transition kernel $\overline{\pi
}^\theta$. Finally, note that the stationarity of this process follows
from a straightforward generalization of Lemma \ref{kiliyir}.
\end{pf}
\begin{lemma}\label{yariuzay}
$\mathbb{Q}_\xi\ll\mathbb{P}$ on $\mathcal{B}_N^+(e_1)$ for every $N\geq
0$; see (\ref{subsigsub}).
\end{lemma}
\begin{pf}
For any $N\geq0$, take an $f\in L^\infty(\mathbb{P})$ such that $f$ is
nonnegative and $\mathcal{B}_N^M$-measurable for some $M\geq0$. Let
$J:=N+M+1$. It follows from (\ref{muyucananan}) and the Schwarz
inequality that
%
%
\begin{eqnarray}\label{seler}\hspace*{15pt}
&&
E_o[\tau_1\exp\{\langle\theta,X_{\tau_1}\rangle- \Lambda_a(\theta)\tau
_1\},\beta=\infty]\int f\,\dd\mathbb{Q}_\xi\nonumber\\
&&\qquad =E_o \Biggl[ \Biggl(\sum_{j=\tau_N}^{\tau_{N+1}-1}f(T_{X_j}\omega
) \Biggr)\exp\{\langle\theta,X_{\tau_J}\rangle- \Lambda_a(\theta)\tau
_J\},\beta=\infty\Biggr]\nonumber\\
&&\qquad \leq\sum_{k=1}^\infty E_o \biggl[\tau_{N+1}=k, \biggl(\sum_{|x|\leq
k}f(T_x\omega) \biggr)\exp\{\langle\theta,X_{\tau_J}\rangle- \Lambda
_a(\theta)\tau_J\},\beta=\infty\biggr]\nonumber\\
&&\qquad =\sum_{k=1}^\infty\sum_{l=1}^k\mathbb{E} \biggl\{ \biggl(\sum
_{|x|\leq k}f(T_x\omega) \biggr)\nonumber\\
&&\qquad\quad\hspace*{43.8pt}{}\times E_o^\omega[\tau_{N+1}=k=H_l,\exp\{
\langle\theta,X_{\tau_J}\rangle- \Lambda_a(\theta)\tau_J\},\beta=\infty
] \biggr\}\nonumber\\
&&\qquad \leq\sum_{k=1}^\infty\sum_{l=1}^k(2k+1)^d\|f\|_{L^2(\mathbb
{P})}\\
&&\qquad\quad\hspace*{27.8pt}{}\times\mathbb{E} \bigl\{E_o^\omega[\tau_{N+1}=k=H_l,\nonumber\\
&&\qquad\quad\hspace*{53.4pt}\hspace*{18.5pt}\exp\{\langle
\theta,X_{\tau_J}\rangle- \Lambda_a(\theta)\tau_J\},\beta=\infty
]^2 \bigr\}^{1/2}.\nonumber
\end{eqnarray}

There exist constants $C_N'<\infty$ and $a_4>0$ such that, for every
$k\geq1$ and $l\in\{1,\ldots,k\}$,
%
%
\begin{eqnarray}\hspace*{28pt}
&&\mathbb{E} \bigl\{E_o^\omega[\tau_{N+1}=k=H_l,\exp\{\langle\theta
,X_{\tau_J}\rangle- \Lambda_a(\theta)\tau_J\},\beta=\infty
]^2 \bigr\}\nonumber\\
\label{seri}
&&\qquad \leq\mathbb{E} \bigl\{E_o^\omega[\tau_{N+1}=k=H_l,\exp\{
\langle\theta,X_{\tau_{N+1}}\rangle- \Lambda_a(\theta)\tau_{N+1}\}
,\beta=\infty]^2 \bigr\}\\
&&\qquad\quad{} \times\Bigl(\inf_{z\in\mathbb{V}_d}P_{o,z}(l_1=0)
\Bigr)^{-1}\nonumber\\
&&\qquad\quad\hspace*{0pt}{}\times\mathbb{E} \bigl\{E_o^\omega[\exp\{\langle\theta,X_{\tau
_M}\rangle- \Lambda_a(\theta)\tau_M\},\beta=\infty]^2 \bigr\}
\nonumber\\
\label{denbil}
&&\qquad \leq C_N'\ee^{-a_4k}.
\end{eqnarray}
\textit{Explanation}: the first term in (\ref{seri}) is bounded from above by
$C_N\ee^{-a_4k}$ for some $C_N<\infty$ and $a_4>0$. The second
term in (\ref{seri}) is finite by Lemma \ref{nonintersectlemma}. Using
the technique in the proof of Lemma \ref{eskiL2lemma}, the third term
in (\ref{seri}) can be shown to be bounded from above by a constant
that is independent of $M$. We leave the details to the reader.

Plugging (\ref{denbil}) into (\ref{seler}), we see that
\[
\int f\dd\mathbb{Q}_\xi\leq C_N''\|f\|_{L^2(\mathbb{P})}
\]
for some finite constant $C_N''$ that is independent of $M$. Since the
functions we have considered are dense in $L^2(\Omega,\mathcal
{B}_N^+(e_1),\mathbb{P})$, it follows from the Riesz representation
theorem that
%
%
\begin{equation}\label{riesz}
\frac{\dd\mathbb{Q}_\xi}{\dd\mathbb{P}}
\bigg|_{\mathcal{B}_N^+(e_1)}\in L^2(\mathbb{P}).
\end{equation}
\upqed\end{pf}

Combining all of the results in this section, we get the following proof.
\begin{pf*}{Proof of Theorem \ref{varminvar}}
Recall that $\mathbb{Q}_\xi\in M_1(\Omega)$ denotes the marginal of
$\hat{\mu}_\xi^\infty\in M_1(\Omega\times U^\mathbb{N})$ which, in
turn, is defined in (\ref{muyucananan}). We have seen that:
\begin{longlist}
\item$\mathbb{Q}_\xi$ is $\hat{\pi}^\theta$-invariant (see Lemma
\ref{kiliyir});
\item$\mathbb{Q}_\xi\ll\mathbb{P}$ on $\mathcal{B}_n^+(e_1)$ for
every $n\geq0$; see Lemma \ref{yariuzay}.
\end{longlist}
It follows from Lemma \ref{Kozlov} (stated below) that $\mathbb{Q}_\xi$
is the unique element of $M_1(\Omega)$ that satisfies these two
properties. We have also proven that $\hat{\mu}_\xi^\infty$ induces an
$\Omega$-valued stationary Markov process with marginal $\mathbb{Q}_\xi
$ and transition kernel $\overline{\pi}^\theta$; see Lemma \ref
{markovkopek}. These results imply part (a) of Theorem \ref{varminvar}.

Note that part (b) of Theorem \ref{varminvar} is a special case of
Theorem \ref{cizme} since $\mathcal{A}_{\eq}\subset\mathcal{A}_a$.
\end{pf*}

In the proof above, we used the following generalization of a classical
homogenization result which is originally due to Kozlov \cite{Kozlov}.
%
\begin{lemma}[(Rassoul-Agha \cite{FirasLLN03})]\label{Kozlov}
Given any $\mathbb{Q}\in M_1(\Omega)$ and any environment kernel $\hat
{\pi}$, define a measure $\mu\in M_1(\Omega\times U)$ by setting
\[
\dd\mu(\cdot,z):=\dd\mathbb{Q}(\cdot)\hat{\pi}(\cdot,z)
\]
for each $z\in U$. Recall Definition \ref{omurtukendi}. If $\mu\in
M_1''(\Omega\times U)$, then the following hold:
%
\begin{enumerate}[(a)]
\item[(a)] the measures $\mathbb{P}$ and $\mathbb{Q}$ are, in fact,
mutually absolutely continuous on $\mathcal{B}_n^+(\hat{v})$ for every
$n\geq0$;
\item[(b)] the environment Markov chain with kernel $\hat{\pi}$ and
initial distribution $\mathbb{Q}$ is stationary and ergodic;
\item[(c)] $\mathbb{Q}$ is the unique $\hat{\pi}$-invariant probability
measure on $\Omega$ that satisfies $\mathbb{Q}\ll\mathbb{P}$ on
$\mathcal{B}_n^+(\hat{v})$ for every $n\geq0$;
\item[(d)] the following LLN is satisfied: $P_o^{\hat{\pi}} (\lim
_{n\rightarrow\infty}\frac{X_n}{n}=\int\sum_{z\in U}\hat{\pi}(\omega
,z)z \,\dd\mathbb{Q} ) = 1$.
\end{enumerate}
\end{lemma}

%
\section{Proof of our results on quenched large deviations}\label{kangiuc}

\subsection{Equality of the quenched and the averaged minimizers}

We start this section by stating the quenched version of Theorem
\ref{averagedconditioningsh}.\vadjust{\goodbreak}
\begin{theorem}\label{stroong}
Assume (\ref{kondisin}). For every $\xi\in\mathcal{A}_{\eq}$,
$\varepsilon
>0$, $N,M,K\geq0$ and $f\dvtx\Omega\times U^{\mathbb{N}}\rightarrow
\mathbb
{R}$ bounded such that $f(\cdot,(z_i)_{i\geq1})$ is independent of
$(z_i)_{i>K}$ and is $\mathcal{B}_N^M$-measurable, the following holds:
%
%
\begin{eqnarray}\label{saglikcan}
&&\limsup_{\delta'\to0}\limsup_{n\rightarrow\infty}\frac{1}{n}\log
P_o^\omega\biggl(\biggl|\int f\,\dd\nu_{n,X}^\infty-\int
f\,\dd\hat{\mu}_\xi^\infty\biggr|>\varepsilon \Big|\nonumber\\
&&\qquad\hspace*{119pt}\hspace*{41.6pt} \biggl|\frac
{X_n}{n}-\xi\biggr|\leq\delta' \biggr)<0\\
\eqntext{\mbox{for $\mathbb
{P}$-a.e. $\omega$.}}
\end{eqnarray}
\end{theorem}
\begin{pf}
This is Theorem 3 of \cite{YilmazSpaceTime}, except that \cite
{YilmazSpaceTime} is concerned with space--time RWRE. However, the
result is obtained directly from Theorem \ref{averagedconditioningsh}
by a standard application of the Borel--Cantelli lemma and Chebyshev's
inequality. In other words, the proof in \cite{YilmazSpaceTime} makes
no use of the space--time assumption. [The only notational difference is
that, in the space--time case, $\Lambda_a(\theta)$ is equal to $\log\phi
(\theta)$ for some explicit function $\phi(\cdot)$, but this does not
play any role in the proof.]
\end{pf}

Now, we are ready to give the following proof.
\begin{pf*}{Proof of Theorem \ref{rosminros}}
Take any $\xi\in\mathcal{A}_{\eq}$. Recall (\ref{Axizz}). If an $\hat
{\alpha}\in A_\xi^\infty$ is not equal to $\hat{\mu}_\xi^\infty$, then
\[
\biggl|\int f\,\dd\hat{\alpha}-\int f\,\dd\hat{\mu}_\xi
^\infty\biggr|>\varepsilon
\]
for some $\varepsilon>0$, $N,M,K\geq0$ and $f\dvtx\Omega\times
U^{\mathbb
{N}}\rightarrow\mathbb{R}$ bounded such that $f(\cdot,(z_i)_{i\geq1})$
is independent of $(z_i)_{i>K}$ and $\mathcal{B}_N^M$-measurable.

For every $\delta'>0$ and $\mathbb{P}$-a.e. $\omega$, (the lower bound
of) the quenched level-3 LDP (i.e., Theorem \ref{level3LDP}) implies that
\[
-I_{q,3}(\hat{\alpha})\leq\liminf_{n\rightarrow\infty}\frac{1}{n}\log
P_o^\omega\biggl(\biggl|\int f\,\dd\nu_{n,X}^\infty-\int
f\,\dd\hat{\mu}_\xi^\infty\biggr|>\varepsilon, \biggl|\frac{X_n}{n}-\xi
\biggr|<\delta' \biggr).
\]
On the other hand,
\[
\lim_{\delta'\rightarrow0}\lim_{n\rightarrow\infty}\frac{1}{n}\log
P_o^\omega\biggl(\biggl|\frac{X_n}{n}-\xi\biggr|\leq\delta' \biggr)=-I_q(\xi),
\]
by the quenched level-1 LDP (i.e., Theorem \ref{qLDPgeneric}). Therefore,
\begin{eqnarray*}
&&
-I_{q,3}(\hat{\alpha})+I_q(\xi)\\
&&\qquad\leq\limsup_{\delta'\to0}\limsup
_{n\rightarrow\infty}\frac{1}{n}\log P_o^\omega\biggl(\biggl|\int
f\,\dd\nu_{n,X}^\infty-\int f\,\dd\hat{\mu}_\xi^\infty
\biggr|>\varepsilon \Big| \biggl|\frac{X_n}{n}-\xi\biggr|\leq\delta' \biggr)\\
&&\qquad<0,
\end{eqnarray*}
by (\ref{saglikcan}). In words, $\hat{\alpha}$ is not a minimizer of
(\ref{level1ratezz}). However, since $I_{q,3}$ is lower semicontinuous
and $A_\xi^\infty$ is compact, there is a minimizer. We conclude that
$\hat{\mu}_\xi^\infty$ is the unique minimizer of (\ref{level1ratezz}).
\end{pf*}
%

\subsection{Modifying Rosenbluth's variational formula}

\begin{lemma}
Assume (\ref{ellipticity}). Recall (\ref{halitakca}) and Definition
\ref{omurtukendi}. For every $\mu\in M_1''(\Omega\times U)$,
%
%
\begin{equation}\label{binsolc}
\mathfrak{I}_q^{**}(\mu)\leq H(\mu).
\end{equation}
\end{lemma}
\begin{pf}
Fix a sequence of test functions, denoted by $(f_i)_{i\geq1}$, that
separate $M_1(\Omega\times U)$. For every $i\geq1$ and $z\in U$, assume
that $f_i(\cdot,z)\dvtx\Omega\to\mathbb{R}$ is measurable with respect to
$\sigma(\omega_x\dvtx|x|\leq i)$. Take any $\mu\in M_1''(\Omega\times U)$.
For every $N\geq1$,
\[
G_{\mu,N}:= \biggl\{\nu\in M_1(\Omega\times U)\dvtx \biggl|\int f_i\,\dd\nu-\int
f_i\,\dd\mu\biggr|<\frac{1}{N}\ \forall i\in\{1,\ldots
,N\} \biggr\}
\]
is an open set. Recall $\hat{v}\in\mathcal{S}^{d-1}$ and the
environment kernel $\hat{\pi}$ corresponding to $\mu$; see Definition
\ref{omurtukendi}. Let $\mathbb{Q}:=(\mu)^1$ so that $\dd\mu
(\cdot,z)=\dd\mathbb{Q}(\cdot)\hat{\pi}(\cdot,z)$ for each $z\in
U$. For every $n\geq1$, introduce a new measure $R_{o,n}^{\hat{\pi
},\omega}$ by setting
\[
\dd R_{o,n}^{\hat{\pi},\omega}:=\frac{\one_{\nu_{n,X}\in G_{\mu
,N},\beta>n}}{P_o^{\hat{\pi},\omega}(\nu_{n,X}\in G_{\mu,N},\beta>n)}
\,\dd P_o^{\hat{\pi},\omega},
\]
where $\beta=\beta(\hat{v}):=\inf\{i\geq0\dvtx \langle X_i,\hat
{v}\rangle<\langle X_o,\hat{v}\rangle\}$.
With this notation, for $\mathbb{Q}$-a.e. $\omega$,
%
%
\begin{eqnarray}\label{evrenak}
&&\log P_o^\omega(\nu_{n,X}\in G_{\mu,N},\beta>n)\nonumber\\
&&\qquad =\log E_o^{\hat{\pi},\omega} \biggl[\nu_{n,X}\in G_{\mu,N},\beta
>n, \frac{\dd P_o^\omega}{\dd P_o^{\hat{\pi},\omega}}
\biggr]\nonumber\\
&&\qquad =\log P_o^{\hat{\pi},\omega}(\nu_{n,X}\in G_{\mu,N},\beta
>n)+\log\int\frac{\dd P_o^\omega}{\dd P_o^{\hat{\pi},\omega
}}\,\dd R_{o,n}^{\hat{\pi},\omega}\nonumber\\
&&\qquad \geq\log P_o^{\hat{\pi},\omega}(\nu_{n,X}\in G_{\mu,N},\beta>n)-
\int\log\frac{\dd P_o^{\hat{\pi},\omega}}{\dd P_o^\omega
}\,\dd R_{o,n}^{\hat{\pi},\omega}\nonumber\\
&&\qquad =\log P_o^{\hat{\pi},\omega}(\nu_{n,X}\in G_{\mu,N},\beta
>n)\\
&&\qquad\quad{} -\frac{1}{P_o^{\hat{\pi},\omega}(\nu_{n,X}\in G_{\mu
,N},\beta>n)} \nonumber\\
&&\qquad\quad\hspace*{10.8pt}{}\times E_o^{\hat{\pi},\omega} \biggl[\nu_{n,X}\in G_{\mu,N},\beta
>n, \log\frac{\dd P_o^{\hat{\pi},\omega}}{\dd P_o^\omega
} \biggr],\nonumber
\end{eqnarray}
by a change of measure and Jensen's inequality.

It follows from Lemma \ref{Kozlov} and the ergodic theorem that
\[
\mathbb{Q}\otimes P_o^{\hat{\pi},\omega}(\nu_{n,X}\in G_{\mu,N}
\mbox{ for sufficiently large $n$})=1
\]
and
\[
\mathbb{Q}\otimes P_o^{\hat{\pi},\omega} \biggl(\lim_{n\to\infty}\frac
{1}{n}\log\frac{\dd P_o^{\hat{\pi},\omega}}{\dd P_o^\omega
}(X_1,\ldots,X_n)=H(\mu) \biggr)=1.
\]
Hence, for $\mathbb{Q}$-a.e. $\omega$,
%
%
\begin{eqnarray}\label{kabusbit}
&\lim_{n\rightarrow\infty}P_o^{\hat{\pi},\omega}(\nu_{n,X}\in G_{\mu
,N},\beta>n)=P_o^{\hat{\pi},\omega}(\beta=\infty)
\end{eqnarray}
and
\begin{eqnarray}
\label{kabusbiter}
&&\limsup_{n\rightarrow\infty}E_o^{\hat{\pi},\omega} \biggl[\nu_{n,X}\in
G_{\mu,N},\beta>n, \frac{1}{n}\log\frac{\dd P_o^{\hat{\pi},\omega
}}{\dd P_o^\omega} \biggr]\nonumber\\[-8pt]\\[-8pt]
&&\qquad\leq P_o^{\hat{\pi},\omega}(\beta=\infty
)H(\mu).\nonumber
\end{eqnarray}
Here, (\ref{kabusbiter}) follows from Fatou's lemma since
\[
\frac{1}{n}\log\frac{\dd P_o^{\hat{\pi},\omega}}{\dd P_o^\omega
}(X_1,\ldots,X_n)=\frac{1}{n}\sum_{i=0}^{n-1}\log\frac{\hat
{\pi}(T_{X_i}\omega,Z_{i+1})}{\pi(X_i,X_{i+1})}\leq-\log\delta,
\]
by uniform ellipticity; see (\ref{ellipticity}).

It follows from parts (b) and (c) of Definition \ref{omurtukendi} that
$P_o^{\hat{\pi},\omega}(\beta=\infty)>0$ for $\mathbb{P}$-a.e. $\omega
$. Since $P_o^{\hat{\pi},\omega}(\beta=\infty)$ is $\mathcal{B}_o^+(\hat
{v})$-measurable, part (d) of Definition \ref{omurtukendi} implies that
%
%
\begin{equation}\label{kaptis}
P_o^{\hat{\pi},\omega}(\beta=\infty)>0 \qquad\mbox{for $\mathbb{Q}$-a.e.
$\omega$.}
\end{equation}

Combining (\ref{evrenak}), (\ref{kabusbit}), (\ref{kabusbiter}) and
(\ref{kaptis}), we see that
%
%
\begin{equation}\label{binbok}
\liminf_{n\rightarrow\infty}\frac{1}{n}\log P_o^\omega(\nu_{n,X}\in
G_{\mu,N},\beta>n)\geq-H(\mu)
\end{equation}
for $\mathbb{Q}$-a.e. $\omega$. However, since $P_o^\omega(\nu
_{n,X}\in G_{\mu,N},\beta>n)$ is $\mathcal{B}_N^+(\hat{v})$-measurable
for every $n\geq1$, Lemma \ref{Kozlov} implies that (\ref{binbok})
holds for $\mathbb{P}$-a.e. $\omega$ as well. Therefore,
\[
\liminf_{n\rightarrow\infty}\frac{1}{n}\log P_o^\omega(\nu_{n,X}\in
G_{\mu,N})\geq-H(\mu) \qquad\mbox{for $\mathbb{P}$-a.e. $\omega$.}
\]

For every $N\geq1$ and $\mathbb{P}$-a.e. $\omega$,
\[
\limsup_{n\rightarrow\infty}\frac{1}{n}\log P_o^\omega(\nu_{n,X}\in
G_{\mu,N})\leq-\inf_{\nu\in\overline{G_{\mu,N}}}\mathfrak{I}_q^{**}(\nu)
\]
by the quenched level-2 LDP, that is, Theorem \ref{level2LDP}. Hence,
\[
\inf_{\nu\in\overline{G_{\mu,N}}}\mathfrak{I}_q^{**}(\nu)\leq H(\mu).
\]
Sending $N$ to infinity implies (\ref{binsolc}) since $\mathfrak
{I}_q^{**}(\cdot)$ is lower semicontinuous and $(f_i)_{i\geq1}$
separates $M_1(\Omega\times U)$.
\end{pf}
\begin{pf*}{Proof of Theorem \ref{gobenk}}
Fix $\xi\neq0$. For any $\mu\in A_\xi\cap M_1'(\Omega\times U)$, let
$\hat{\pi}$ be the environment kernel given by $\hat{\pi}(\cdot
,z):=\frac{\dd\mu(\cdot,z)}{\dd(\mu)^1(\cdot)}$ for each
$z\in U$. It is shown in \cite{Kozlov} that
\[
P_o^{\hat{\pi}} \biggl(\lim_{n\to\infty}\frac{X_n}{n}=\xi\biggr)=1.
\]
Hence, $\mu\in M_1''(\Omega\times U)$. [For part (c) of Definition \ref
{omurtukendi}, take any $\hat{v}\in S^{d-1}$ such that $\langle\xi,\hat
{v}\rangle>0$.] In other words,
%
%
\begin{equation}\label{yavuzz}
A_\xi\cap M_1'(\Omega\times U)\subset A_\xi\cap M_1''(\Omega\times U).
\end{equation}
It follows from (\ref{level1rate}), (\ref{binsolc}), (\ref{yavuzz}) and
(\ref{guzelyaziss}) that
\begin{eqnarray*}
I_q(\xi)&=&\inf_{\mu\in A_\xi}\mathfrak{I}_q^{**}(\mu)\leq\inf\{\mathfrak
{I}_q^{**}(\mu)\dvtx \mu\in A_\xi\cap M_1''(\Omega\times U)\}\\
&\leq&\inf\{H(\mu)\dvtx \mu\in A_\xi\cap M_1''(\Omega\times U)\}\\
&\leq&\inf\{H(\mu)\dvtx \mu\in A_\xi\cap M_1'(\Omega\times U)\}\\
&=&I_q(\xi).
\end{eqnarray*}
In particular, $I_q(\xi)=\inf\{H(\mu)\dvtx \mu\in A_\xi\cap
M_1''(\Omega
\times U)\}$. 
\end{pf*}
\begin{pf*}{Proof of Theorem \ref{tezbitti}}
Fix $\xi\in\mathcal{A}_{\eq}$. Then, $\xi\neq0$. Indeed, by
differentiating both sides of (\ref{yalebbim}) with respect to $\theta$
at $\theta=\nabla I_a(\xi)$, we see that
\[
\langle\xi,e_1\rangle=\langle\nabla\Lambda_a(\theta),e_1\rangle=\frac
{E_o[\langle X_{\tau_1},e_1\rangle\exp\{\langle\theta,X_{\tau_1}\rangle
-\Lambda_a(\theta)\tau_1\}|\beta=\infty]}{E_o[\tau_1\exp\{\langle\theta
,X_{\tau_1}\rangle-\Lambda_a(\theta)\tau_1\}|\beta=\infty]}>0.
\]
%

Recall that $\hat{\mu}_\xi\in M_1(\Omega\times U)$ is the marginal of
$\hat{\mu}_\xi^\infty\in M_1(\Omega\times U^\mathbb{N})$ and $\dd\hat
{\mu}_\xi(\cdot,z)=\dd\mathbb{Q}_\xi(\cdot)\hat{\pi
}^\theta(\cdot,z)$ for each $z\in U$. It is clear from Theorem \ref
{varminvar} that $\hat{\mu}_\xi\in A_\xi\cap M_1''(\Omega\times U)$.
Observe that
%
%
\begin{eqnarray}\label{cinkobir}\quad
H(\hat{\mu}_\xi)&=&\int\sum_{z\in U}\hat{\pi}^\theta(\omega,z)\log\frac
{\hat{\pi}^\theta(\omega,z)}{\pi(0,z)}\,\dd\mathbb{Q}_\xi(\omega
)\\
\label{cinkoiki}
&=&\int\sum_{z\in U}\hat{\pi}^\theta(\omega,z) \biggl(\langle\theta
,z\rangle-\Lambda_a(\theta)+\log\frac{h(\theta,T_z\omega)}{h(\theta
,\omega)} \biggr)\,\dd\mathbb{Q}_\xi(\omega)\\
\label{cinkouc}
&=&\langle\theta,\xi\rangle-\Lambda_a(\theta)+\int\sum_{z\in U}\hat{\pi
}^\theta(\omega,z)\log\frac{h(\theta,T_z\omega)}{h(\theta,\omega
)}\,\dd\mathbb{Q}_\xi(\omega)\\
\label{cinkodort}
&=&\langle\theta,\xi\rangle-\Lambda_a(\theta)\\
\label{cinkobes}
&=&I_q(\xi).
\end{eqnarray}
\textit{Explanation}: (\ref{cinkobir}), (\ref{cinkoiki}) and (\ref{cinkouc})
follow from (\ref{halitakca}), (\ref{azdinyn}) and (\ref{ximu}),
respectively. Since $\mathbb{Q}_\xi$ is $\hat{\pi}^\theta$-invariant by
Lemma \ref{kiliyir}, it is easy to see that the integral in (\ref
{cinkouc}) is zero. Finally, (\ref{cinkodort}) is equal to (\ref
{cinkobes}) because $\xi=\nabla\Lambda_a(\theta)$ and $I_q(\xi)=I_a(\xi
)$. 

Thus far, we have shown that $H(\hat{\mu}_\xi)=I_q(\xi)$. Now, take any
$\nu\in A_\xi\cap M_1''(\Omega\times U)$. If $\nu\neq\hat{\mu}_\xi$,
then
\[
I_q(\xi)<\mathfrak{I}_q^{**}(\nu)\leq H(\nu),
\]
by Corollary \ref{duzbac} and (\ref{binsolc}).
We conclude that $\hat{\mu}_\xi$ is the unique minimizer of (\ref
{guzelyazisss}).
\end{pf*}
%


%

%
\printaddresses

\end{document}